\documentclass[11pt]{article}
\usepackage{psfrag,
graphicx, harvard,
epsf,mathrsfs,
euscript,amsfonts,amsmath,latexsym,amssymb,mathrsfs
}
\usepackage[active]{srcltx}
\usepackage{color}

\usepackage{float}
\usepackage{multirow}

\newtheorem{theorem}{Theorem}

\newtheorem{remark}{Remark}

\newtheorem{corollary}{Corollary}

\newcommand{\cA}{{\cal A}}

\newcommand{\bR}{{\mathbb R}}

\newcommand{\sF}{{\mathscr F}}

\newcommand{\sR}{{\mathscr R}}

\newcommand{\sT}{{\mathscr T}}

\newcommand{\sX}{{\mathscr X}}
\newcommand{\sY}{{\mathscr Y}}
\newcommand{\sZ}{{\mathscr Z}}

\newcommand{\rP}{\mathrm{P}}
\newcommand{\rE}{\mathrm{E}}

\newcommand{\rd}{\mathrm{d}}

\newcommand{\pr} {\par \noindent{\bf Proof\,:~}}
\newcommand{\epr}{\hfill\hbox{\hskip 4pt
                \vrule width 5pt height 6pt depth 1.5pt}\vspace{0.5cm}\par}
\begin{document}
\title{A Unified
Approach
for Solving  Sequential Selection Problems
}
\author{
{\sf Alexander Goldenshluger}\thanks{Department of Statistics, University of Haifa, Haifa 31905, Israel.  e-mail:
{\em goldensh@stat.haifa.ac.il}.}
\\ University of Haifa
\and
\and {\sf Yaakov Malinovsky}
\thanks{Department of Mathematics and Statistics,
University of Maryland, Baltimore County, MD 21250, USA.
e-mail: {\em yaakovm@umbc.edu}}
\\ University of Maryland
\\
Baltimore County
\and
\and {\sf Assaf Zeevi}\thanks{Graduate School of Business, Columbia University, New York, NY 10027, USA.  e-mail: {\em assaf@gsb.columbia.edu}}
\\ Columbia University
}
\date{}
\maketitle

\begin{abstract}
In this paper we develop a unified
approach
for solving
a wide class of sequential selection problems.  This class includes, but is not limited to,  selection problems with no--information,
rank--dependent rewards, and considers both  fixed as well as random problem horizons.
The proposed framework
is based on a reduction of the original selection problem to one of
{\it optimal stopping} for a sequence of judiciously constructed independent random variables. We demonstrate that our
approach
allows
exact and efficient
computation  of optimal policies and various performance metrics thereof for a variety of sequential
selection problems, several of which have not been solved to date.

\end{abstract}
\vspace*{1em}
\par
\vspace*{1em}

\section{Introduction}
\label{se:I}
In
sequential selection problems a decision maker examines  a sequence of observations  which appear in random order over some horizon.  Each observation  can be either accepted or rejected, and  these decisions are irrevocable. The objective is to select an element in this sequence  to optimize a given criterion.  A classical example is the  so-called {\em secretary problem} in which the objective is to maximize the probability of selecting the element of the sequence that ranks highest.
The existing literature
contains numerous settings and formulations of such problems, see, e.g.,
\citeasnoun{GiMo}, \citeasnoun{freeman},
\citeasnoun{gnedin-book}, \citeasnoun{ferguson},
\citeasnoun{samuels} and \citeasnoun{F2008}; to make more concrete connections we defer further references to the subsequent section where we formulate the class of problems more precisely.
\par
Sequential selection problems are typically
solved using the principles of dynamic programming, relying heavily on structure that is problem-specific, and focusing on theoretical properties of the optimal solution; cf. \citeasnoun{GiMo}, \citeasnoun{gnedin-book} and
\citeasnoun{F2008}.
Consequently, it has become increasingly difficult to discern  commonalities among the multitude of problem variants and their solutions. Moreover, the resulting optimal policies are often viewed as difficult to implement, and focus is  placed on deriving
sub--optimal policies and various asymptotic approximations; see, e.g.,
\citeasnoun{mucci-a},  \citeasnoun{FrSa},
\citeasnoun{krieger-ester}, and \citeasnoun{Arlotto}, among many others.
\par
In this paper we demonstrate that a wide
class of such problems can be solved optimally and in a unified manner.
 This class includes, but is not limited to, sequential selection problems
with {\em no--information}, {\em rank--dependent rewards}
and allows for fixed or random horizons.
The proposed solution methodology covers both problems that have been worked out in the literature, albeit in an instance-specific manner, as well as several problems whose solution to the best of our knowledge is not known to date.
 We refer
to Section~\ref{sec:selection-problems} for details.
%
The unified framework we develop is based on the fact that various
sequential selection problems can be
reduced, via a conditioning argument,  to a problem of optimal stopping for a sequence of independent random variables
that are constructed in a special way.
The latter is an instance  of a more general class of problems, referred to as {\it sequential stochastic assignments}, first formulated and solved by \citeasnoun{DLR} (some extensions are given in  \citeasnoun{albright}).  The main idea of the proposed framework
was briefly
sketched  in
\citeasnoun[Section~4]{GZ}; in this paper it is fully fleshed and adapted to the range of problems alluded to above.
\par
The approach we take is operational, insofar as it supports exact and efficient
computation of the  optimal policies and corresponding optimal values, as well as various other performance metrics.
In the  words of
\citeasnoun{robbins70}, we
 ``put the problem on a computer.''
Optimal stopping rules that result from our
approach
belong to the class of   memoryless threshold policies and hence have a relatively simple structure.  In particular,
the proposed reduction
constructs a new sequence of {\it independent} random variables, and the optimal  
rule is to stop  the first time instant when the current ``observation'' exceeds 
a given threshold. The threshold 
computation is predicated on the structure of  the policy in sequential stochastic assignment problems
\`a  la \citeasnoun{DLR} and \citeasnoun{albright} (as part of the so pursued unification, these problems are also extended in the present paper to the case of a random time horizon).  The structure of the optimal stopping rule we derive allows us to explicitly compute probabilistic characteristics  and various performance metrics of the stopping time, which, outside of special cases, are completely absent from the literature.

\par
The  rest of the paper is structured as follows.
Section \ref{sec:selection-problems}
discusses sequential selection problems. In this section we formulate
two general no--information problems with rank--dependent reward corresponding to  fixed and random horizon
[Problems~(A1) and~(A2) respectively]. We  also present various specific problem instances, Problems~(P1)--(P12), that are
covered by the proposed unified framework. Section \ref{sec:SAP} describes the class of stochastic sequential selection problems:
we consider the standard formulation, Problem~(AP1), first introduced and solved by \citeasnoun{DLR},
and a formulation with random horizon, Problem~(AP2). These problems
are central to our solution approach. Section~\ref{sec:B}  presents the auxiliary stopping problem, Problem~(B),
and explains its solution via the mapping to a stochastic assignment problem. It then explains the details of the reduction and the structure of the algorithm that implements our proposed stopping rule.
Section~\ref{sec:examples} presents the implementation of said algorithm to
Problems~(P1)--(P12)
surveyed in
Secton~\ref{sec:selection-problems}.
We close with a few concluding remarks in Section \ref{sec:conc}.

%

\section{Sequential selection problems}\label{sec:selection-problems}
\par
Let us introduce some notation and terminology.
Let $X_1, X_2, \ldots$ be an infinite  sequence of independent identically distributed continuous
random variables defined on a probability space $(\Omega, \sF, \rP)$.
Let $R_t$  be the relative rank of $X_t$ and  $A_{t,n}$ be the absolute rank of $X_t$ among the first $n$ observations (which we also refer to as the {\it problem horizon}):
\begin{equation}\label{eq:ranks}
 R_t:=\sum_{j=1}^t {\bf 1}(X_t\leq X_j),\;\;
A_{t,n}:=\sum_{j=1}^n {\bf 1}(X_t\leq X_j),\;\;\;t=1,\ldots,n.
\end{equation}
Note that
with this  notation the largest observation  has the absolute rank one, and
$R_t=A_{t,t}$ for any $t$.
Let $\sR_t:=\sigma(R_1,\ldots,R_t)$ and $\sX_t:=\sigma(X_1,\ldots,X_t)$ denote
the $\sigma$--fields generated by $R_1,\ldots,R_t$ and $X_1,\ldots,X_t$,  respectively;
$\sR=(\sR_t, 1\leq t\leq n)$ and $\sX=(\sX_t, 1\leq t\leq n)$ are the corresponding filtrations.
In general, the class of all stopping times of a filtration $\sY=(\sY_t, 1 \leq t\leq n)$
will be denoted $\sT(\sY)$; i.e.,
$\tau\in \sT(\sY)$ if $\{\tau=t\}\in \sY_t$ for all $1\leq t\leq n$.
 \par
 Sequential selection problems are classified according to
 the information available to the decision maker and the structure of the reward function.
The settings in which only relative ranks $\{R_t\}$ are observed
 are usually referred to as \hbox{\em no--information problems}, whereas
 {\em full information} refers to the case when random variables
 $\{X_t\}$  can be observed,
 and their distribution is  known. 
 In addition,
 the total number of available observations $n$
 can be either fixed or random with given distribution.
These cases are referred to  as problems with fixed  and random horizon,
 respectively.

 \subsection{Problems with fixed horizon}\label{sec:fixed-horizon}
In this paper we mainly consider 
selection problems with {\em no--information} and
{\em rank--dependent reward}.
 The prototypical  sequential
 selection problem with
 fixed horizon,
 no--information and rank--dependent reward is formulated as follows; see, e.g.,
 \citeasnoun{gnedin-krengel}.
 \begin{quote}
 {\sc Problem~(A1)}.
 Let $n$ be a fixed positive integer, and let $q:\{1,2,\ldots, n\}\to \bR$ be a reward function.
 The average reward of a stopping rule $\tau\in \sT(\sR)$ is
 \[
 V_n(q; \tau):=\rE q\big(A_{\tau,n}\big).
 \]
The objective is
 to find  the  rule
 $\tau_*\in \sT(\sR)$ satisfying
 \[
 V^*_n(q) := \max_{\tau\in \sT(\sR)} V_n(q; \tau)= \rE q\big(A_{{\tau_*,n}}\big)
 \]
 and to compute the optimal value
$V^*_n(q)$.
\end{quote}
 Depending on the reward function $q$ we distinguish among the following
 types of sequential selection problems with fixed horizon.
%
 \paragraph{Best--choice problems.}
 The settings in which the reward function is an indicator
 are usually referred to  as  {\em best--choice stopping problems}. Of special note are the following.
 \par\medskip
 (P1).~{\em Classical secretary problem.} This
 problem
setting corresponds to the case $q(a)=q_{\rm csp}(a):={\bf 1}\{a=1\}$.
 Here we want to maximize the probability  $\rP\{A_{\tau,n}=1\}$
of selecting the best
alternative over all
stopping times $\tau$ from $\sT(\sR)$. It is well known that
the optimal policy will pass on
approximately the first $n/e$ observations and select the first
subsequent to that which is superior than all previous ones, if such an observation exists; otherwise the last element in the sequence  is selected.
The limiting optimal value is
$\lim_{n\to\infty}V^*_n(q_{\rm csp})=1/e$ \cite{lindley,Dyn,GiMo}.
Ferguson~(1989) reviews the problem history  and discusses how different
assumptions about this problem evolved over time.
\par\medskip
(P2).~{\em Selecting  one of the $k$ best values.} The problem is usually referred to as
{\em the Gusein--Zade stopping problem}  \cite{GuZa,FrSa}.
Here  $q(a)=q_{\rm gz}^{(k)}(a):={\bf 1}\{a\leq k\}$, and the problem is
to maximize $\rP\{A_{\tau,n} \leq k\}$ with respect to
$\tau\in \sT(\sR)$.
 The optimal policy was characterized in
\citeasnoun{GuZa}. It is determined by $k$ natural numbers $1\leq \pi_1\leq \cdots\leq \pi_k$
and proceeds as follows: pass the first $\pi_1-1$ observations
and among the subsequent $\pi_1, \pi_1+1,\ldots,\pi_2-1$ observations
choose
the first observation with relative rank one;
if it does not exists then among the set of observations $\pi_2, \pi_2+1,\ldots,\pi_3-1$  choose
the one of  relative rank two,
etc.
\citeasnoun{GuZa} presented dynamic programming algorithm to determine the numbers $\pi_1,\ldots,\pi_k$ and
value of $V^*_n(q_{\rm gz}^{(k)})$.
He also studied the limiting behavior of the numbers $\pi_1,\ldots,\pi_k$ as the problem horizon grows large, and showed that
$\lim_{n\to\infty}V^*_n(q_{\rm gz}^{(2)})\approx 0.574$.
Exact results
for the case $k = 3$ are given in \citeasnoun{QuLaw} and for general $k$ in \citeasnoun{Woryna2017}.

Based on
general asymptotic results of \citeasnoun{mucci-a}, \citeasnoun{FrSa} computed numerically
$\lim_{n\to\infty} V^*_n\big(q_{\rm gz}^{(k)}\big)$
for a range of different values of $k$. The recent paper
\citeasnoun{DiLaRi} studies some approximate policies.
\par\medskip
(P3).~{\em Selecting  the $k$th best alternative.}
In this problem $q(a)=q_{\rm pd}^{(k)}(a):=
{\bf 1}\{a=k\}$, i.e.
we want to maximize
the probability of selecting the $k$th best candidate. The problem was explicitly solved for $k=2$
by \citeasnoun{Szajowski1982}, \citeasnoun{Rose} and \citeasnoun{Vanderbei2012}; the last paper coined the name
the {\em postdoc  problem} for this setting. An optimal policy for $k=2$ is to reject first
$\lceil n/2\rceil$ observations and then  select the one which is the second best
relative to this previous observation set, if it exits; otherwise the last element in the sequence  is selected. The optimal value is
\mbox{$V^*_n(q_{\rm pd}^{(2)})=(n+1)/4n$} if $n$ is odd and \mbox{$V^*_n(q_{\rm pd}^{(2)})=n/4(n-1)$} if $n$ is even.
An optimal stopping rule for the case $k=3$ and some results
on the optimal value were reported recently in \citeasnoun{Yao}.
We are not aware of results
on the optimal policy and exact computation of the optimal values
for general $n$ and $k$.
Recently approximate policies were developed in \citeasnoun{Bruss-2016}.
The problem of selecting the median value $k=(n+1)/2$,
where $n$ is odd, was
considered in \citeasnoun{Rose-2}. It is shown there that
$\lim_{n\to\infty} V^*_n(q_{\rm pd}^{((n+1)/2)})=0$.
\paragraph{Expected rank type problems.}
To this category we  attribute problems
with  reward $q$ which is not an indicator function.
\par\medskip
(P4).~{\em Minimization of the expected rank.}
In this problem the goal is  to
minimize $\rE A_{\tau, n}$ with respect to $\tau\in\sT(\sR)$.
If we put $q(a)=q_{\rm er}(a):=-a$ then
\begin{equation}\label{eq:chow-rob}
 \min_{\tau\in \sT(\sR)} \rE A_{\tau,n} =
 - \max_{\tau\in \sT(\sR)} \rE \,q_{\rm er}(A_{\tau,n}).
\end{equation}
This problem was discussed heuristically by
\citeasnoun{lindley} and solved by \citeasnoun{chow}.
It was shown there that $\lim_{n\to\infty} \min_{\tau\in\sT(\sR)} \rE A_{\tau,n}
=\prod_{j=1}^\infty(1+\frac{2}{j})^{1/(j+1)}\approx 3.8695$.
The corresponding  optimal stopping rule is given by backward
induction relations.
A simple suboptimal stopping rule which is close
to the optimal one was proposed in  \citeasnoun{krieger-ester}.
\par\medskip
(P5).~{\em Minimization of the expected squared rank.}
Based on \citeasnoun{chow},
\citeasnoun{Robbins-91}
developed the optimal policy and computed
the asymptotic optimal value in  the problem of minimization of
$\rE [A_{\tau,n}(A_{\tau,n}+1)\cdots (A_{\tau,n}+k-1)]$ with respect to $\tau\in \sT(\sR)$.
In particular, he showed that for the optimal stopping rule $\tau_*$
\[
 \lim_{n\to\infty}\rE[A_{\tau_*,n}(A_{\tau_*,n}+1)\cdots (A_{\tau_*,n}+k-1)] =
 k! \bigg\{\prod_{j=1}^\infty
\bigg(1+ \frac{k+1}{j}\bigg)^{1/(k+j)}\bigg\}^k.
 \]
\citeasnoun{Robbins-91}
also discussed the problem of minimization of $\rE A^2_{\tau,n}$ over
$\tau\in \sT(\sR)$ and
mentioned that the optimal stopping rule and optimal value are unknown.
As we will demonstrate  below,
optimal policies for any problem of this type can be easily derived, and the corresponding
optimal values are straightforwardly calculated for any fixed $n$.

\subsection{Problems with random horizon}
In Problem~(A1) and specific problem instances of Section~\ref{sec:fixed-horizon}
the horizon  $n$
is fixed beforehand, and optimal policies depend critically on this assumption.
However, in practical situations $n$ may be  unknown. This fact motivates
settings in which the horizon is assumed to be
a random variable.
\par
If the horizon is random then
the selection may not have been made by the time 
the observation process terminates.
In order to  take this possibility into account, we introduce minor modifications in the definitions of
the absolute and relative ranks in (\ref{eq:ranks}).
By convention we put $A_{t,k}=0$ for $t>k$, and if
$N$ is a positive random variable representing the problem horizon and taking values in
$\{1,\ldots, N_{\max}\}$
($N_{\max}$ can be infinite) then
on the event $\{N=k\}$, $k=1, \ldots, N_{\max}$,  we set
\begin{equation}\label{eq:barR}
                        \bar{R}_t:=\left\{ \begin{array}{ll}
                         R_{t}, & t=1, \ldots, k,\\
                         0& t=k+1, \ldots, N_{\max}.
                        \end{array}
\right.
\end{equation}
Furthermore, $\bar{\sR}_t:=\sigma(\bar{R}_1, \dots, \bar{R}_t)$ denotes
the $\sigma$--field induced by $(\bar{R}_1, \ldots, \bar{R}_t)$, and
$\bar{\sR}:=\{\bar{\sR}_t, 1\leq t\leq N_{\max}\}$ is the corresponding filtration. We refer to the sequence
$\{\bar{R}_t, 1\leq t\leq N_{\max}\}$ as the sequence of observed relative ranks.
\par
The general
 selection problem with
  random horizon,
 no--information and rank--dependent reward is formulated as follows [see
 \citeasnoun{PS1972} and \citeasnoun{Irle80}].
\begin{quote}
 {\sc Problem~(A2)}.
Let  $N$ be a positive integer  random
variable
with distribution $\gamma=\{\gamma_k\}$, $\gamma_k=\rP(N=k)$, $k=1,2,\ldots, N_{\max}$, where
$N_{\max}$ may be infinite.
Assume that $N$ is independent of the sequence $\{X_t, t\geq 1\}$.
Let
$q:\{1,\ldots, N_{\max}\}\cup \{0\} \to [0, \infty)$ be a reward function,
and by convention $q(0)=0$.
Let $Q_t:=q(A_{t,k})$, $t\in\{1,\ldots, N_{\max}\}$ on the event $\{N=k\}$.
The performance
of a stopping rule $\tau\in \sT(\bar{\sR})$ is measured~by
$V_\gamma(q;\tau):= \rE Q_{\tau}$.
The objective is to  find  the stopping rule
 $\tau_*\in \sT(\bar{\sR})$ such that
 \[
 V_\gamma^*(q) := \max_{\tau\in \sT(\bar{\sR})} \rE Q_{\tau}=
 \max_{\tau\in \sT(\bar{\sR})} V_\gamma(q; \tau)=
  V_\gamma(q;\tau_*)
  \]
  and to compute  the optimal value $V^*_\gamma(q)$.
\end{quote}
The introduced  model
assigns  fictitious zero value to the observed relative rank $\bar{R}_t$ if
the selection has not been made by the end of the problem
horizon, i.e., if  $t>N$. By assumption $q(0)=0$
the reward for not selecting an observation by time $N$ is also set to zero, though
other possibilities can be considered
for this value.
\par\medskip
In principle, all
problems (P1)--(P5)  discussed above can be  formulated and solved under the assumption that the observation horizon is random.
Below we discuss the following  three problem instances.
\par\medskip
(P6).~{\em Classical secretary problem with random horizon.}
The classical secretary problem with random horizon $N$
corresponds to Problem~(A2) with $q(a)={\bf 1}\{a=1\}$; it
was studied in \citeasnoun{PS1972}.
In  Problem (P1) where $n$ is fixed,
the stopping region is an interval of the form
$\{k_n, \ldots, n\}$ for some integer $k_n$.
In contrast to (P1),
\citeasnoun{PS1972} show that
for general distributions of $N$
the optimal policy can involve ``islands,'' i.e., the stopping region can
be a union of several disjoint intervals (``islands'').
The paper derives some sufficient conditions under which
the stopping region is a single interval and presents specific examples
satisfying these conditions.
In particular, it is shown that in the  case
of the uniform distribution  on $\{1,\ldots, N_{\max}\}$, i.e.,
$\gamma_k=1/N_{\max}$, $k=1,\ldots, N_{\max}$,
the stopping region  is of the form $\{k_{N_{\max}},\ldots, N_{\max}\}$ with
$k_{N_{\max}}/N_{\max} \to
 2e^{-2}$,
 $V^*_{\gamma}(q_{\rm csp})\to
 2e^{-2}$ as $N_{\max}\to\infty$.
The characterization of
optimal policies for general distributions of $N$ is not available in the existing
literature.
\par\medskip
(P7).~{\em Selecting  one of the $k$ best values over a random horizon.}
This is a version of the Gusein--Zade stopping problem, Problem~(P2), with
random horizon. Recall that here the reward function is $q_{\rm gz}^{(k)}(a)={\bf 1}\{a\leq k\}$.
To the best of our knowledge, this setting has been studied only for  $k=2$ and uniform
distribution of $N$, i.e., $\gamma_k=1/N_{\max}$, $k\in \{1, \ldots, N_{\max}\}$; see \citeasnoun{KT2003}.
The cited paper derives the optimal policy and demonstrates that
it is qualitatively the same as in the setting with  fixed horizon.
\citeasnoun{KT2003} study asymptotics of thresholds $\pi_1$ and $\pi_2$,
and compute numerically  the problem optimal value for a range of $N_{\max}$'s; in particular,
$\lim_{N_{\max}\to\infty} V^*_\gamma(q_{\rm gz}^{(2)})\approx 0.4038$.
Below we show how this problem can be stated and solved for general $k$ and arbitrary distribution
of $N$ within our proposed unified framework.
\par\medskip
(P8). {\em Minimization of the expected rank over a random horizon.}
Consider a variant of Problem~(P4) under the assumption that the horizon
is a random variable $N$ with known distribution.
In this setting the loss (the negative reward) for stopping at time
$t$ is the absolute rank $A_{t, N}$ on the event $\{N\geq t\}$; otherwise, the absolute rank of the last available observation
$A_{N, N}=R_N$ is received.
We want to minimize the expected loss over all stopping rules
$\tau\in \sT(\bar{\sR})$.
This problem has been considered in \citeasnoun{Gianini-Pettitt}.
In particular, it was shown there that if $N$ is uniformly distributed over $\{1,\ldots, N_{\max}\}$ then
the expected loss tends to infinity as $N_{\max}\to \infty$.
On the other hand, for distributions which are more
``concentrated'' around $N_{\max}$, the optimal value coincides
asymptotically with the one for Problem~(P4).
Below we demonstrate that this problem can be naturally  formulated and solved
for general distributions of $N$ using our  proposed unified framework;
 the details are
given in Section~\ref{sec:Pettitt}.

\subsection{Multiple choice problems}
The proposed framework is also applicable for some multiple choice problems
both with fixed and random horizons. Below we review  two settings
with fixed horizon.
\par\medskip
(P9). {\em Maximizing the  probability of
selecting the best observation  with $k$ choices.}
Assume that one can  make $k$ selections, and the reward function
equals one if the best observation belongs to the selected subset and zero otherwise.
Formally, the problem
is to maximize the probability $\rP(\cup_{j=1}^k \{A_{\tau_j ,n}=1\})$ over
stopping times $\tau_1<\cdots<\tau_k$ from $\sT(\sR)$. This problem
has been considered in  \citeasnoun{GiMo} who gave numerical results for up to
$k=8$; see also \citeasnoun{Haggstrom} for theoretical results for $k=2$.
\par\medskip
(P10). {\em Minimization of the expected average rank.}
Assume that $k$ choices are possible, and the goal is to minimize the expected average
rank of the selected subset. Formally, the problem is  to minimize
$\frac{1}{k} \rE \sum_{j=1}^k A_{\tau_j,n}$ over stopping times $\tau_1<\cdots<\tau_k$
of $\sT(\sR)$. For related results we refer to
\citeasnoun{Megiddo},
\citeasnoun{kep-jap}, \citeasnoun{kep-aap1} and
\citeasnoun{Nikolaev-Sofronov}.

\subsection{Miscellaneous problems}
The proposed framework extends beyond problems with rank--dependent rewards and no--information.
The next two problem instances demonstrate such extensions.
\par\medskip
(P11). {\em Moser's problem with random horizon.}
Let $\{X_t, t\geq 1\}$ be a sequence of independent identically distributed
random variables with distribution $G$ and expectation $\mu$.
Let $N$ be a positive integer--valued random variable representing the problem horizon.
We observe sequentially $X_1, X_2, \ldots$ and the reward for stopping at time $t$ is the  value
of the observed random variable $X_t$; if the stopping does not occur by problem horizon $N$, then the reward
is the last observed observation $X_N$. Formally, we want to maximize
\[
 \rE\big[X_{\tau}{\bf 1}\{\tau\leq N\}+X_{N}{\bf 1}\{\tau>N\}\big],
\]
with respect to all stopping times $\tau$ of the filtration associated with the observed values.
The formulation with fixed $N=n$
and uniformly distributed $X_t$'s on $[0,1]$ corresponds to the classical problem of \citeasnoun{Moser}.
\par\medskip
(P12).
{\em Bruss'  Odds--Theorem.}
\citeasnoun{Bruss} considered the following optimal stopping problem.
Let $Z_1, \ldots, Z_n$ be independent Benoulli random variables with success probabilities
$p_1,\ldots, p_n$ respectively. We observe $Z_1, Z_2,\ldots$ sequentially and  want to
stop at the time of the last success, i.e., the problem is to
find a stopping time
$\tau \in \sT(\sZ)$ so as the probability
$\rP(Z_\tau=1, Z_{\tau+1}=Z_{\tau+2}=\cdots=Z_n=0)$ is maximized.
Odds--Theorem \cite[Theorem~1]{Bruss} states that it is optimal to stop
at the first time instance $t$ such  that
\[
Z_t=1\;\;\;\hbox{and}\;\;\;
 t\geq t_*:=\sup\bigg\{1, \,\sup\Big\{ k=1,\ldots, n: \sum_{j=k}^n \frac{p_j}{q_j} \geq 1\Big\}\bigg\},
\]
with $q_j:=1-p_j$ and $\sup\{\emptyset\}=-\infty$.
This statement has been used in various settings for finding optimal stopping policies.
For example, it provides 
shortest
self--contained solution to the classical
secretary problem \cite{Bruss}. For some extensions  to multiple stopping
problems see \citeasnoun{Matsui} and references therein. We also refer to
the recent work
\citeasnoun{Bruss19} where  further relevant references can be found.
In what follows we will demonstrate that Bruss' Odds--Theorem can be
derived using the proposed framework.
\section{Sequential stochastic assignment problems}\label{sec:SAP}
The unified framework we propose leverages the sequential assignment model toward
the solution of the problems presented in Section~\ref{sec:selection-problems}.
In this section we consider two formulations of the stochastic sequential
assignment problem: the first is the classical formulation introduced by
\citeasnoun{DLR}, while the second one is an extension for random
horizon.

\subsection{Sequential assignment problem with fixed horizon}
The formulation below follows the terminology used by \citeasnoun{DLR}. Suppose
that $n$ jobs arrive sequentially in time, referring henceforth to the latter as the problem horizon.  The $t$th job, $1\leq t\leq n$, is identified with
a random variable $Y_t$ which is observed. The jobs must be assigned to $n$ persons
which have known
``values'' $p_1,\ldots, p_n$.  Exactly one job should be assigned to
each person, and after the assignment the person becomes unavailable for
the next jobs. If the $t$th job is assigned to the $j$th person then a reward of $p_jY_t$
is obtained.
The goal is to maximize the expected total reward.
\par
Formally, assume that $Y_1,\ldots, Y_n$ are integrable independent random variables defined on
probability space $(\Omega, \sF, \rP)$, and let $F_t$ be the distribution function of $Y_t$ for each $t$.
Let  $\sY_t$ denote the $\sigma$--field generated by $(Y_1,\ldots, Y_t)$:
$\sY_t=\sigma(Y_1,\ldots,Y_t)$,
$1\leq t\leq n$.
Suppose that $\pi=(\pi_1,\ldots, \pi_n)$ is a permutation of
$\{1,\ldots,n\}$ defined on $(\Omega,\sF)$. We say that $\pi$
is an {\em assignment policy} (or simply {\em policy}) if
$\{\pi_t=j\}\in \sY_t$
for every  $1\leq j\leq n$ and  $1\leq t\leq n$.
That is, $\pi$ is a policy if it is non--anticipating
relative to the filtration $\sY=\{\sY_t, 1\leq t\leq n\}$ so that
$t$th job is assigned on the basis of information in $\sY_t$.
Denote by $\Pi(\sY)$ the set of all policies
associated with the filtration $\sY=\{\sY_t, 1\leq t\leq n\}$.

\par\medskip
Now consider the following sequential assignment problem.
\begin{quote}
{\sc Problem (AP1).}
Given a vector $p=(p_1,\ldots,p_n)$, with $p_1\leq p_2\leq\cdots\leq p_n$,  we want
to maximize {\em  the total expected reward}
$S_{n}(\pi) :=\rE\sum_{t=1}^n p_{\pi_t} Y_t$
with respect to $\pi\in \Pi(\sY)$.
The policy $\pi^*$ is called {\em optimal} if
$S_n(\pi^*)  = \sup_{\pi\in \Pi(\sY)} S_n(\pi)$.
\end{quote}
In the sequel the following
representation will be useful
\[
\sum_{t=1}^n p_{\pi_t} Y_t = \sum_{t=1}^n \sum_{j=1}^n p_j Y_t {\bf 1}\{\pi_t=j\} =
\sum_{j=1}^n p_j Y_{\nu_j};
\]
here the random variables $\nu_j\in \{1,\ldots,n\}$, $j=1,\ldots, n$ are given by the
one-to-one correspondence
$\{\nu_j=t\}=\{\pi_t=j\}$, $1\leq t\leq n$, $1\leq j\leq n$.
In words, $\nu_j$ denotes
the index of the job to which the $j$th person
is assigned.
\par
The structure of the optimal policy is given by the following statement.
%
\begin{theorem}[\citeasnoun{DLR}; \citeasnoun{albright}]
\label{th:derman}
Consider Problem~(AP1) with horizon $n$.
There exist real numbers $\{a_{j,n}\}_{j=0}^n$,
\[
 -\infty\equiv a_{0,n}\leq a_{1,n}\leq \cdots\leq a_{n-1,n}\leq a_{n,n}\equiv \infty
\]
such that on the first step,
when random variable $Y_1$ distributed $F_1$ is observed,
the optimal policy is $\pi^*_1=\sum_{j=1}^n j {\bf 1}\{Y_1\in (a_{j-1,n},a_{j,n}]\}$.
 The numbers $\{a_{j,n}\}_{j=1}^n$ do not depend on $p_1,\ldots,p_n$
and are determined by the following recursive relationship
\[
 a_{j,n+1}=\int_{a_{j-1,n}}^{a_{j,n}} z \rd F_1(z) + a_{j-1,n} F_1(a_{j-1,n}) + a_{j,n}[1-F_1(a_{j,n})],\;\;
j=1,\ldots,n,
\]
where $-\infty\cdot 0$ and $\infty\cdot 0$ are defined to be $0$. At the end of the first stage the assigned $p$ is removed from the feasible set and the process repeats with the next observation, where the above calculation is then performed relative to the distribution $F_2$ and real numbers
$
 -\infty\equiv a_{0,n-1}\leq a_{1,n-1}\leq \cdots\leq a_{n-2,n-1}\leq a_{n-1,n-1}\equiv \infty
$ are determined and so on.
Moreover, $a_{j,n+1}=\rE Y_{\nu_j}$, $\forall 1\leq j\leq n$, i.e.,
$a_{j,n+1}$ is the expected value of the job which is assigned to the $j$th person,
and $\sum_{j=1}^n p_ja_{j, n+1}$ is the optimal value of the problem.
\end{theorem}

\begin{remark}
In order to determine an optimal policy we calculate inductively a triangular array
$\{a_{j,t}\}_{j=1}^{t-1}$ for $t=2,\ldots, n+1$, where   $F_{n-t+2}$ is used in order to compute $\{a_{j,t}\}_{j=1}^{t-1}$.
In implementation the  optimal policy uses numbers $a_{1,n},\,a_{2,n},\,a_{n-1,n}$  in order to identify one value from $p_1,\ldots, p_n$
which will multiply  $Y_1$. Then, this value of $p$ is excluded from $n$ values, and numbers
$a_{1,n-1},\,a_{2,n-1},\,a_{n-2,n-1}$ are used for determination
of the next value of $p$ from $n-1$ remaining values; this value will multiply $Y_2$, and so on. At the last step the number
$a_{1,2}$ is to assign one of the two remaining values of $p$ to  $Y_{n-1}$. Finally, the last
remaining value of $p$ will be assigned to  $Y_n$.
\end{remark}

\subsection{Stochastic sequential assignment problems with random horizon}

In practical situations the horizon, or number of available jobs,  $n$ is often unknown. Under these circumstances
 the optimal policy of \citeasnoun{DLR} is not applicable.
This fact provides motivation for the setting  with random number of jobs.
The
sequential assignment problem with random horizon was formulated and solved  by
\citeasnoun{sakaguchi} who derived the optimal policy using dynamic programming principles.
More recently, \citeasnoun{Jacobson} also
considered the sequential assignment problem  with a random horizon.
They show that the optimal solution to the problem with random horizon can be derived
from the solution to an  auxiliary  assignment problem with dependent job sizes.
Below we demonstrate  that the problem with random horizon is in
fact equivalent to  a certain version
of the sequential assignment problem with fixed horizon and independent job sizes.

\par\medskip
The stochastic sequential assignment problem with random horizon is stated as follows.
\begin{quote}
{\sc Problem (AP2).}  Let $N$ be a positive integer-valued  random variable with distribution
$\gamma=\{\gamma_k\}$, $\gamma_k=\rP(N=k)$, $k=1,\ldots, N_{\max}$, where
$N_{\max}$ can be infinite. Let $Y_1, Y_2, \ldots $ be an infinite sequence of integrable independent
random variables with distributions $F_1, F_2, \ldots$
such that $\rP(Y_t=0)=0$ for all $t$.
Assume that $N$ is independent of $\{Y_t, t\geq 1\}$.
Let $\bar{Y}_1, \bar{Y}_2, \ldots$ be the sequence of random variable defined as follows: if $N=k$, $k\in \{1, \ldots, N_{\max}\}$ then
\begin{equation}\label{eq:Y-bar}
 \bar{Y}_t = \left\{\begin{array}{ll}
                       Y_t, & t\leq k,\\
                       0, & t>k,
                      \end{array}
\right.
\;\;\;t=1, 2, \ldots, N_{\max}.
\end{equation}
Let $\bar{\sY}_t:=\sigma(\bar{Y}_1, \ldots, \bar{Y}_t)$ be the $\sigma$--field induced
by $(\bar{Y}_1, \ldots, \bar{Y}_t)$, and
$\bar{\sY}=\{\bar{Y}_t, 1\leq t\leq N_{\max}\}$ be the corresponding filtration.
Given real numbers
$p_1\leq \ldots\leq p_{N_{\max}}$ the objective is to maximize the expected total reward
$S_\gamma(\pi)= \rE\sum_{t=1}^N p_{\pi_t} Y_t$
over all policies $\pi\in \Pi(\bar{\sY})$.
\end{quote}
\par
\begin{remark} \mbox{}
\begin{itemize}
\item[{\rm (i)}]
The probability model  of Problem~(AP2) postulates that  the decision maker observes
 vector $(\bar{Y}_1, \ldots, \bar{Y}_{N_{\max}})$
that is  generated as follows.
Given random variable $N$ and
a  sequence  $\{Y_t, t\geq 1\}$, independent of $N$,
the decision maker is presented with the $N_{\max}$--vector
$(Y_1, \ldots, Y_k, 0, \ldots, 0)$ on the event $\{N=k\}$, $k\in\{1, \ldots, N_{\max}\}$.
Thus, the distribution of $(\bar{Y}_1, \ldots, \bar{Y}_{N_{\max}})$
 is the mixture
of distributions of vectors
\[
 (Y_1, 0,\ldots, 0), \;(Y_1, Y_2, 0,\ldots, 0), \cdots, (Y_1, Y_2, \ldots, Y_{N_{\max}})
\]
with respective weights $\gamma_1$, $\gamma_2, \ldots, \gamma_{N_{\max}}$.
\item[{\rm (ii)}]
The definition of the
sequence $\{\bar{Y}_t, t\geq 1\}$ and condition
$\rP(Y_t=0)=0$ for all $t$ imply that the first observed
zero value of $\bar{Y}_t$ designates termination of the assignment process.
In particular,  $\bar{Y}_t=0$ implies that $\bar{Y}_s=0$ for all $s\geq t$.
\end{itemize}
\end{remark}
\par\medskip
In the following statement we show that Problem~(AP2) is equivalent
to a version of Problem~(AP1), the  standard sequential assignment problem
with fixed horizon and independent job sizes.
\begin{theorem}\label{th:random}
The optimal value in Problem~(AP2) coincides with
the optimal value in Problem~(AP1) associated with fixed horizon $n=N_{\max}$ and
independent job sizes $Y_t\sum_{k=t}^{N_{\max}}\gamma_k$.
The optimal policy in Problem~(AP2) follows the one in   Problem~(AP1) with fixed horizon $n=N_{\max}$ and
independent job sizes $Y_t\sum_{k=t}^{N_{\max}}\gamma_k$ until the first zero value of $\bar{Y}_t$ is observed;
this indicates termination of the assignment process.
\end{theorem}
\pr
With the introduced notation  for any $\pi\in \Pi(\bar{\sY})$
\begin{equation}\label{eq:S}
S_\gamma(\pi)= \rE \sum_{t=1}^N p_{\pi_t} Y_t =  \rE \sum_{t=1}^{N_{\max}} p_{\pi_t} \bar{Y}_t
=  \sum_{t=1}^{N_{\max}} \rE \big[p_{\pi_t} Y_t {\bf 1}\{N\geq t\}\big].
\end{equation}
It follows from (\ref{eq:S})  that the expected total reward
$S_\gamma(\pi)$ is fully determined  by
the values of $p_{\pi_t}$ on
events $\{N\geq t\}$,  $t=1,\ldots, N_{\max}$ only; the value
of $p_{\pi_t}$ on $\{N<t\}$
is irrelevant as the ensuing reward is equal to zero.
Note that
$\pi_t$ is $\bar{\sY}_t$--measurable, i.e., $\pi_t=\pi_t(\bar{Y}_1, \ldots, \bar{Y}_t)$
for any $t=1, \ldots, N_{\max}$.
However,  by definition,  $\bar{Y}_1=Y_1,\ldots, \bar{Y}_t=Y_t$ on the event $\{N\geq t\}$;
hence $\bar{\sY}_t\cap \{N\geq t\}= \sY_t\cap \{N\geq t\}$, and
$\pi_t=\pi_t(Y_1, \ldots, Y_t)$ on $\{N\geq t\}$. This implies that in
(\ref{eq:S}) the decision variable
$\pi_t$ can be taken to be $\sY_t$--measurable.
It follows that
\begin{eqnarray*}
 \rE\big[p_{\pi_t} Y_t {\bf 1}\{N\geq t\}\big] = \rE\Big\{
 \rE \big[p_{\pi_t} Y_t \,{\bf 1}\{N\geq t\} | \sY_t\big]\Big\}=
 \rE \Big\{p_{\pi_t} Y_t \sum_{k=t}^{N_{\max}} \gamma_k\Big\},
\end{eqnarray*}
where
the last equality follows from independence of $N$ and $\sY_t$.
Thus,
\[
 S_\gamma(\pi) = \rE \sum_{t=1}^{N_{\max}} p_{\pi_t} \Big\{Y_t \sum_{k=t}^{N_{\max}} \gamma_k
 \Big\},
\]
which shows that
the optimal value coincides with the one in the assignment problem with fixed
horizon $n=N_{\max}$ and independent job sizes $Y_t\sum_{k=t}^{N_{\max}} \gamma_k$.
As long as the assignment process proceeds,  the optimal policy
follows the one in said problem with fixed horizon $n=N_{\max}$ and independent job sizes
$Y_t\sum_{k=t}^{N_{\max}} \gamma_k$. The first observed
zero value of $\bar{Y}_t$ indicates termination of the assignment
process due to horizon randomness.
\epr
\begin{remark}

To the best of our knowledge, the relation between Problems~(AP2) and~(AP1)
established in Theorem~\ref{th:random} is new.
In fact,
this relationship is implicit in the optimal policy derived in
\citeasnoun{sakaguchi}; however, \citeasnoun{sakaguchi} does not mention this.
In contrast, \citeasnoun{Jacobson}
develop optimal policy by
reduction of the problem to an
auxiliary one with
dependent job sizes.  As
 Theorem~\ref{th:random} shows, this is not necessary:
 the problem with random number of jobs is equivalent to the standard sequential assignment
 problem with independent job sizes, and  it is solved by the standard procedure
 of \citeasnoun{DLR}.
\end{remark}

\begin{remark}
In Theorem~\ref{th:random} we assume that $N_{\max}$ is finite. Under suitable assumptions on
the  weights $\{p_j\}$ and jobs sizes $\{Y_t\}$ one can construct $\epsilon$--optimal policies
 for the problem with infinite $N_{\max}$. However,  we do not pursue this direction here.
\end{remark}

\section{A unified approach
for solving sequential selection problems}\label{sec:B}
\subsection{An auxiliary optimal stopping problem}
Consider the following auxiliary problem of optimal stopping.
\begin{quote}
{\sc Problem~(B)}. Let $Y_1,\ldots, Y_n$ be a sequence of integrable
independent real-valued  random variables
with corresponding distributions
$F_1, \ldots, F_n$.
For a stopping rule $\tau\in \sT(\sY)$ define
$W_n(\tau):=\rE Y_\tau$. The objective is to
find the stopping rule $\tau_*\in \sT(\sY)$ such that
\[
 W_n^*:=\max_{\tau\in \sT(\sY)} \rE Y_\tau= W_n(\tau_*)=\rE Y_{\tau_*}.
\]
\end{quote}
\par
Problem~(B) is a specific case of the stochastic sequential assignment problem
of \citeasnoun{DLR},
and Theorem~\ref{th:derman}
has immediate implications for Problem~(B).
The following statement is a  straightforward consequence of Theorem~\ref{th:derman}.
\begin{corollary}\label{cor:opt-stop}
Consider Problem~(B). Let
$\{b_t, \,t\geq 1\}$ be the sequence of real numbers defined recursively
by
\begin{align}
&b_1=-\infty,\,\,\,b_2=\rE Y_{n},\,\,\,\,\,
\nonumber
\\
&b_{t+1}= \int_{b_t}^\infty z \rd F_{n-t+1}(z) + b_t F_{n-t+1}(b_t),\,\,\,t=2,\dots,n.
\label{eq:b}
\end{align}
Let
\begin{equation}\label{eq:opt-rule}
\tau_*= \min\{1\leq t\leq n: Y_t > b_{n-t+1}\};
\end{equation}
then
\[
W_n^*=\rE Y_{\tau_*}=\max_{\tau\in \sT(\sY)} \rE Y_\tau = b_{n+1}.
\]
\end{corollary}

\pr The integral in (\ref{eq:b}) is finite because the  random variables
$Y_1,\ldots, Y_n$ are integrable.
Consider Problem~(AP1) with
$p=(0,\ldots,0,1)$.
By Theorem~\ref{th:derman}, at step $t$
the optimal policy assigns value $p_n$ to the job $Y_t$ only if
$Y_t > a_{n-t, n-t+1}$, $t=1,\ldots,n$, and
\[
 a_{n-t, n-t+1} = \int_{a_{n-t-1, n-t}}^\infty z \rd F_{t+1}(z) + a_{n-t-1, n-t} F_{t+1}(a_{n-t-1,n-t}).
\]
Setting $b_t:=a_{t-1,t}$, and noting that 
 $b_1=-\infty,\,\,b_2=\int_{-\infty}^\infty z \rd F_{n}(z)$,
we come to the required statement. \epr

%
\subsection{Reduction to the auxiliary stopping problem} \label{sec:reduction}
Problems~(A1) and (A2) of Section~\ref{sec:selection-problems}
can be reduced to the optimal
stopping of a sequence of independent random variables [Problem~(B)].
In order to demonstrate this relationship we use well known properties
of the relative and absolute ranks defined in (\ref{eq:ranks}).
These properties are briefly recalled
in the next paragraph; for details see, e.g., \citeasnoun{gnedin-krengel}.
\par
Let $A_n:=(A_{1,n},\ldots, A_{n,n})$, and
let $\cA_n$ denote then set of all permutations of $\{1,\ldots, n\}$; then
$\rP(A_n =A)= 1/n!$ for all $A\in \cA_n$ and all $n$.
The random variables $\{R_t, t\geq 1\}$ are independent, and
 $\rP(R_t=r) = 1/t$ for all  $r=1,\ldots, t$.
For any $n$ and $t=1,\ldots, n$
\begin{equation}\label{eq:A|R0}
\rP(A_{t,n}=a|R_1=r_1,\ldots,R_t=r_t)=\rP(A_{t,n}=a|R_t=r_t),
\end{equation}
and
\begin{eqnarray}\label{eq:A|R}
&& \rP(A_{t,n}=a|R_t=r) \;=\;
\frac{\binom{a-1}{r-1}\binom{n-a}{t-r}}{\binom{n}{t}},\;\;\;r\leq a\leq n-t+r.
\end{eqnarray}

\par
Now we are in a position to establish a relationship between Problems~(A1) and~(B).
\paragraph{Fixed horizon.}
Let
\begin{eqnarray}\label{eq:u-t}
I_{t,n}(r):=
\sum_{a=r}^{n-t+r} q(a) \frac{\binom{a-1}{r-1}\binom{n-a}{t-r}}{\binom{n}{t}},
\;\;\;\; r=1,\ldots, t.
\end{eqnarray}
It follows from (\ref{eq:A|R})  that
$I_{t,n}(R_t)= \rE\{ q(A_{t,n})\, |\, R_t\}$.
Define
\begin{eqnarray}
\label{eq:Y}
 Y_t := I_{t,n}(R_t),\;\;\;t=1,\ldots, n.
\end{eqnarray}
By independence of the relative ranks,  $\{Y_t\}$ is a sequence of
independent random variables.
\par
The relationship between stopping problems~(A1) and~(B) is given in the next theorem.
\begin{theorem}\label{th:reduction}
The optimal
stopping rule $\tau_*$ solving
Problem~(B) with
random variables $\{Y_t\}$ given in (\ref{eq:u-t})--(\ref{eq:Y})
also solves Problem~(A1):
\[
V_n(q; \tau_*)=\max_{\tau\in \sT(\sR)} \rE q(A_{\tau,n})=
 \max_{\tau\in \sT(\sY)} \rE Y_\tau= W_n(\tau_*).
\]
\end{theorem}
\pr
First we note that for any stopping rule $\tau\in \sT(\sR)$
one has
$\rE q(A_{\tau,n}) = \rE Y_\tau$,  where $Y_t:= \rE[q(A_{t,n})|\sR_t]$.
Indeed,
\begin{eqnarray*}
 \rE q(A_\tau) &=& \sum_{k=1}^n
 \rE q(A_\tau) {\bf 1}\{\tau=k\}=\sum_{k=1}^n
 \rE q(A_k) {\bf 1} \{\tau=k\}
 \\
 &=& \sum_{k=1}^n
 \rE\Big[ {\bf 1}\{\tau=k\} \rE\{q(A_k)|\sR_k\}\Big] = \sum_{k=1}^n \rE [{\bf 1}\{\tau=k\} Y_k] =
 \rE Y_\tau,
\end{eqnarray*}
where we have used the fact that $\{\tau=k\}\in \sR_k$. This implies that
$\max_{\tau\in \sT(\sR)} \rE q(A_{\tau,n})=\max_{\tau\in \sT(\sR)} \rE Y_\tau$.
To prove the theorem it suffices to show only that
\begin{equation}\label{eq:reduction}
\max_{\tau\in \sT(\sR)} \rE Y_\tau=\max_{\tau\in \sT(\sY)}\rE Y_\tau.
\end{equation}
Clearly,
\begin{equation}
\label{eq:C1}
\sY_t\subset \sR_t, \;\;\;\forall 1\leq t\leq n.
\end{equation}
Because $R_1,\ldots,R_n$ are independent random variables, and
$Y_t=I_{t,n}(R_t)$, $\forall t$
we have that for any $s, t\in \{1,\ldots,n\}$ with $s<t$
\begin{equation}
\label{eq:C2}
 \rP\{G_t\,|\,\sY_s\}=\rP\{G_t\,|\,\sR_s\},\;\;\;\forall G_t\in \sY_t.
\end{equation}
The statement \eqref{eq:reduction} follows from \eqref{eq:C1}, \eqref{eq:C2}
and Theorem~5.3 of \citeasnoun{chow-rob-sieg}.
In fact, \eqref{eq:reduction} is a consequence
of the well known fact that randomization does not increase rewards
in stopping problems \cite[Chapter~5]{chow-rob-sieg}. This concludes the proof.
\epr
It follows from Theorem~\ref{th:reduction} that the optimal stopping rule
in Problem~(A1) is given by Corollary~\ref{cor:opt-stop} with random variables $\{Y_t\}$ defined
by (\ref{eq:Y}). To implement  the rule we need to  compute   the distributions $\{F_t\}$ of the random variables
$\{Y_t\}$ and to apply formulas (\ref{eq:b}) and~(\ref{eq:opt-rule}).

\paragraph{Random horizon.} 

Next, we establish a correspondence between Problems~(A2) and~(B).
Let
\begin{equation}\label{eq:tilde-u}
 J_t(r) :=\sum_{k=t}^{N_{\max}} \gamma_k I_{t,k}(r), \;\;\;r=1,\ldots, t,
\end{equation}
where $I_{t,k}(\cdot)$ is given
in (\ref{eq:u-t}), and $\gamma_k=\rP(N=k)$.
Below in the proof  of Theorem~\ref{th:random-reduction}
we show that
\[
J_{t}(r)= \rE \big\{ q(A_{t,N}) {\bf 1}\{N\geq t\} | R_1=r_1, \ldots, R_{t-1}=r_{t-1}, R_t=r\big \}.
\]
Define also
\begin{equation}\label{eq:tilde-Y}
 Y_t:= J_t(R_t)=\sum_{k=t}^{N_{\max}} \gamma_k I_{t,k}(R_t),\;\;\;t=1,\ldots, N_{\max}.
\end{equation}
\begin{theorem}\label{th:random-reduction}\mbox{}
\begin{itemize}
\item[{\rm (i)}] Let $N_{\max}<\infty$; then the  optimal
stopping rule $\tau_*$  solving
Problem~(B) with fixed horizon $N_{\max}$
and random variables $\{Y_t\}$ given in
(\ref{eq:tilde-u})--(\ref{eq:tilde-Y})
provides
the optimal solution to Problem~(A2):
\[
 V_\gamma^*(q)=\max_{\tau\in \sT(\bar{\sR})} V_\gamma (q;\tau)=
 \max_{\tau\in \sT(\sY)} \rE Y_\tau = W_{N_{\max}} (\tau^*).
\]
\item[{\rm (ii)}] Let $N_{\max}=\infty$ and
assume that
\begin{equation}\label{eq:u-finite}
\sup_t \max_{1\leq r\leq t} \;\sum_{k=t}^\infty \gamma_k |I_{t,k}(r)| <\infty.
\end{equation}
Let
$\epsilon>0$ be arbitrary; then there exists
$\tilde{N}_{\max}=\tilde{N}_{\max}(\epsilon)$ such that
for any stopping rule $\tau\in \sT(\bar{\sR})$ one has
\begin{equation}\label{eq:eps-1}
 W_{\tilde{N}_{\max}}(\tau)-\epsilon\leq V_\gamma (q;\tau) \leq
 W_{\tilde{N}_{\max}}(\tau)+
 \epsilon.
\end{equation}
In particular,
the optimal stopping rule $\tau_*$ solving Problem~(B) with fixed
horizon $\tilde{N}_{\max}= \tilde{N}_{\max}(\epsilon)$ and $\{Y_t\}$ given
(\ref{eq:tilde-u})--(\ref{eq:tilde-Y}) is an $\epsilon$--optimal
stopping rule for Problem~(A2):
\begin{equation}\label{eq:eps-2}
 W_{\tilde{N}_{\max}}(\tau_*)-\epsilon \leq V_\gamma^*(q) \leq
 W_{\tilde{N}_{\max}}(\tau_*) +\epsilon
\end{equation}
\end{itemize}
\end{theorem}
\pr
(i). In Problem~(A2)
the reward  for stopping at time $t$ is
\mbox{$Q_t=q(A_{t, N}) {\bf 1}\{N\geq t\}$}, and the objective is
to maximize $\rE Q_\tau$ with respect to stopping times $\tau$ of filtration
$\bar{\sR}$ [see~(\ref{eq:barR})].
First,
we argue that as long as the decision process does not terminate  before time~$t$,
we can restrict ourselves to stopping times $\tau$
adapted to filtration $\sR$. This is a consequence of the fact that
performance $V_\gamma (q; \tau)=\rE Q_\tau$ of any stopping rule $\tau\in \sT(\bar{\sR})$
is fully determined
by its probabilistic properties on the event $\{\tau\leq N\}$ only.
Indeed, write
\[
 Q_\tau= q(A_{\tau, N})  {\bf 1}\{N\geq \tau\}= \sum_{t=1}^{N_{\max}}
 q(A_{t, N}) {\bf 1}\{\tau=t\} {\bf 1}\{N\geq t\}.
\]
The event $\{\tau=t\}$ belongs to $\bar{\sR}_t$, i.e.,
${\bf 1}\{\tau=t\}=:\varphi_t=\varphi_t(\bar{R}_1, \ldots, \bar{R}_t)$ is a measurable function of
$\bar{R}_1, \ldots, \bar{R}_t$.
However, on the event $\{N\geq t\}$, when
the decision
process is
at time $t$, we have
$\bar{R}_1=R_1, \ldots, \bar{R}_t=R_t$ so that in fact
$\varphi_t=\varphi_t(R_1, \ldots, R_t)$. Thus,
in view of the structure of the reward function, at any time instance
$t$ at which the decision is made
we should consider stopping rules adapted to $\sR$ only, i.e., $\tau\in \sT(\sR)$.
This implies by conditioning
\begin{eqnarray}
 \rE Q_\tau &=&  \rE  \sum_{t=1}^{N_{\max}} \rE\big[ q(A_{t, N}) {\bf 1}\{N\geq t\} {\bf 1}\{\tau=t\}\,|\, \sR_t\big]
\nonumber
 \\
 &=&  \rE  \sum_{t=1}^{N_{\max}} {\bf 1}\{\tau=t\} \rE\big[ q(A_{t, N}) {\bf 1}\{N\geq t\} \,|\, \sR_t\big]
\nonumber
 \\
 &=&
 \rE \sum_{t=1}^{N_{\max}} {\bf 1}(\tau=t)  \sum_{k=t}^{N_{\max}}
 \gamma_k \rE\big[ q(A_{t, k})\,|\, \sR_t\big]
 =\rE Y_\tau,
 \label{eq:Q}
\end{eqnarray}
where $Y_t=\sum_{k=t}^{N_{\max}} \gamma_k I_{t,k}(R_t)$, $t=1,\ldots, N_{\max}$ [cf.~(\ref{eq:tilde-Y})].
Here the second equality follows from  $\{\tau=t\}\in \sR_t$
on $\{N\geq t\}$, while the third equality holds by independence of $N$ and~$\{R_t, t\geq 1\}$.
The remainder of the proof proceeds along the lines of
the proof of Theorem~\ref{th:reduction}.
\par
(ii). In view of the proof of~(i) we can restrict ourselves with with the stopping rules $\tau\in \sT(\sR)$.
Let $\tilde{N}_{\max}=\tilde{N}_{\max}(\epsilon)$ be the minimal integer number such that
\begin{equation}\label{eq:condition}
\sup_t \max_{1\leq r\leq t} \,\sum_{k=\tilde{N}_{\max}+1}^{\infty }
\gamma_k |I_{t,k}(r)| \leq \epsilon.
\end{equation}
The existence of $\tilde{N}_{\max}(\epsilon)$ follows from
(\ref{eq:u-finite}).
By (\ref{eq:Q}) and (\ref{eq:condition}),
for any stopping rule \mbox{$\tau\in \sT(\sR)$}  we have
$V_{\gamma}(q;\tau)=\rE \sum_{k=\tau}^\infty \gamma_k I_{\tau, k}(R_\tau)$,
and
\begin{eqnarray*}
\rE \sum_{k=\tau}^{\tilde{N}_{\max}}  \gamma_k I_{\tau ,k}(R_\tau) -\epsilon
\;\leq\;
V_\gamma(q;\tau) \leq
\rE \sum_{k=\tau}^{\tilde{N}_{\max}}  \gamma_k I_{\tau ,k}(R_\tau) + \epsilon.
\end{eqnarray*}
This implies (\ref{eq:eps-1}).
In order to prove (\ref{eq:eps-2}) we note that if $\tilde{\tau}$
is the optimal stopping rule in Problem~(A2) then by (\ref{eq:eps-1})
and definition of $\tau_*$
\[
 V_\gamma(q; \tilde{\tau})= V_\gamma^*(q)\leq W_{\tilde{N}_{\max}}(\tilde{\tau})
 +\epsilon \leq W_{\tilde{N}_{\max}}(\tau_*)+\epsilon,
\]
which proves the upper bound in (\ref{eq:eps-2}).  On the other hand,
in view of (\ref{eq:eps-1})
\[
 V^*_\gamma(q)=V_\gamma(q; \tilde{\tau}) \geq
 V_\gamma(q; \tau_*) \geq W_{\tilde{N}_{\max}}(\tau_*)-\epsilon.
\]
This concludes the proof.
\epr

\begin{remark}
 Condition~(\ref{eq:u-finite}) imposes restrictions on the tail of the distribution of $N$.
 It can be easily
 verified in any concrete setting; for details see Section~\ref{sec:examples}.
\end{remark}

\begin{remark}
 Theorems~\ref{th:reduction} and~\ref{th:random-reduction}  imply that
solution of
 Problems~(A1) and~(A2) can be obtained by solving
 Problem~(B) with  a suitably defined horizon and random variables
$\{Y_t\}$ given
by (\ref{eq:u-t})--(\ref{eq:Y}) and
(\ref{eq:tilde-u})-(\ref{eq:tilde-Y})
respectively. The latter problem is solved by the recursive procedure
given in Corollary~\ref{cor:opt-stop}.
\end{remark}
\subsection{Specification of the optimal stopping rule for Problems~(A1) and~(A2)}
Now, using Theorems~\ref{th:reduction} and~\ref{th:random-reduction},
we specialize the result of Corollary~\ref{cor:opt-stop}
for solution of Problems~(A1) and~(A2).
For this purpose we require the following notation:
\[
 \nu:=\left\{ \begin{array}{ll}
             n, & {\rm Problem~(A1)},\\
             N_{\max}\hbox{ or } \tilde{N}_{\max}, & {\rm Problem~(A2)},
            \end{array} \right.
\;\;\;
U_t(r):=\left\{ \begin{array}{ll}
             I_{t, n}(r), & {\rm Problem~(A1)},\\
             J_t(r),  & {\rm Problem~(A2)}.
            \end{array} \right.
\]
\par
Note that in Problem~(A2)
we put $\nu=N_{\max}$ for distributions with the finite right endpoint
$N_{\max}<\infty$; otherwise
$\nu=\tilde{N}_{\max}$,
where $\tilde{N}_{\max}$ is defined in the proof of Theorem~\ref{th:random-reduction}.
With this notation Problem~(B) is associated with independent random variables
$Y_t=U_t(R_t)$ for $t=1,\ldots, \nu$.
\par
Let $y_{t}(1), \ldots,y_{t}(\ell_t)$ denote distinct points of the set
$\{U_{t}(1),\ldots,U_{t}(t)\}$, $t=1,\ldots, \nu$.
The distribution of the random variable $Y_t$ is supported
on  the set $\{y_{t}(1),\ldots,y_{t}(\ell_t)\}$ and given~by
\begin{eqnarray}
\label{eq:distr-Y}
 f_t(j)&:=&\rP\{Y_t= y_{t}(j)\}= \frac{1}{t}\sum_{r=1}^t
{\bf 1}\big\{ U_{t}(r)=y_{t}(j)\big\},\;\;\;\;j=1,\ldots,\ell_t,
\\
 F_t(z)&=&\sum_{j=1}^{\ell_t} f_t(j) {\bf 1}\{y_t(j)\leq z\},\;\;\;z\in \bR.
\label{eq:distr-Y1}
\end{eqnarray}
The following statement is an immediate consequence of Corollary~\ref{cor:opt-stop}
and formulas (\ref{eq:distr-Y})--(\ref{eq:distr-Y1}).
\begin{corollary}\label{cor:main}
 Let $\tau_*= \min\{1\leq t\leq \nu: Y_t > b_{\nu-t+1}\}$,
 where the sequence $\{b_t\}$ is given~by
\begin{eqnarray}
&& b_1=-\infty,\,\,\,b_2=\sum_{j=1}^{\ell_{\nu}} y_{\nu}(j) f_{\nu}(j),\,\,\,\,
\label{eq:rec-aa}
\\
&& b_{t+1}=\sum_{j=1}^{\ell_{\nu-t+1}} \big[\,b_t \vee y_{\nu-t+1}(j)\,\big]\,f_{\nu-t+1}(j),\;\;\;
t=2,\dots,\nu.
\label{eq:rec-b}
\end{eqnarray}
Then
\[
 \rE Y_{\tau_*}=\sup_{\tau\in \sT(\sR)} \rE Y_\tau=b_{\nu+1}.
\]
\end{corollary}
\pr
In view of \eqref{eq:Y} and (\ref{eq:tilde-Y}) , $Y_1,\ldots,Y_\nu$ are independent random variables; therefore Corollary~\ref{cor:opt-stop}
is applicable.
We have
\begin{eqnarray*}
\int_{b_t}^\infty z \rd F_{\nu-t+1}(z)&=&
\sum_{j=1}^{\ell_{\nu-t+1}} y_{\nu-t+1}(j){\bf 1}\{y_{\nu-t+1}(j)>b_t\} f_{\nu-t+1}(j),
\\
b_t F_{\nu-t+1}(b_t) &=& b_t \sum_{j=1}^{\ell_{\nu-t+1}}
f_{\nu-t+1}(j) {\bf 1}\{y_{\nu-t+1}(j)\leq b_t\}.
\end{eqnarray*}
Summing up these expressions we come to (\ref{eq:rec-b}).
\epr
\paragraph{Expectation of stopping times.}
As we have already  mentioned, in the considered problems
the optimal stopping rule belongs to the class of  memoryless threshold policies. This facilitates  derivation of the
distributions of the
corresponding stopping times, and calculation of their probabilistic characteristics.
One of the important characteristics is the expected time elapsed before stopping. In  problems with fixed horizon $\nu=n$ it is given by the following formula
\begin{eqnarray}
\label{ex:ExpTimeToStop}
\rE\left(\tau_*\right)&=&\sum_{i=0}^{n-1}\rP(\tau_*>i)=1+\sum_{i=1}^{n-1}\rP(\tau_*>i)
\nonumber
\\
&=& 1+\sum_{i=1}^{n-1}\prod_{t=1}^{i}\rP\left(Y_t\leq b_{n-t+1}\right)=1+\sum_{i=1}^{n-1}\prod_{t=1}^{i}F_{t}\left(b_{n-t+1}\right),
\end{eqnarray}
where $\{F_t\}$ and $\{b_t\}$ are defined in (\ref{eq:distr-Y1}) and
(\ref{eq:rec-aa})--(\ref{eq:rec-b}).
\par
In the problems where the horizon $N$ is random,
the time until stopping is $\tau_*\wedge N$. In this case
\begin{equation}\label{eq:tau-min}
\rE (\tau_*\wedge N)= \rE \tau_* {\bf 1}\{\tau_*\leq N\}+ \rE N {\bf 1}\{\tau_*>N\},
\end{equation}
where
\begin{align}\label{eq:111}
 \rE &[\tau_* {\bf 1}\{\tau_*\leq N\}] = \rE \Big(\tau_* \sum_{k=\tau_*}^{N_{\max}} \gamma_k\Big)
 = \sum_{j=1}^{N_{\max}} j \sum_{k=j}^{N_{\max}} \gamma_k \rP(\tau_*=j)
 \nonumber
 \\
 &= \sum_{k=2}^{N_{\max}} \gamma_k (1-F_1(b_{N_{\max}})) +
 \sum_{k=2}^{N_{\max}} \gamma_k \sum_{j=2}^{k} j (1-F_{j}(b_{N_{\max}-j+1}))\prod_{t=1}^{j-1}
 F_t(b_{N_{\max}-t+1})
\end{align}
and
\begin{align}\label{eq:222}
 \rE[N {\bf 1}(N<\tau_*)] = \sum_{k=1}^{N_{\max}} k \gamma_k \prod_{t=1}^k F_t(b_{N_{\max}-t+1}).
\end{align}

\subsection{Implementation}\label{sec:implementation}
In this section we
present an efficient algorithm implementing the optimal stopping rule described earlier.
In order to
implement  (\ref{eq:rec-aa})--(\ref{eq:rec-b}) we need to
find the sets
$\{y_t(j), j=1,\ldots, \ell_t\}$ in  which
random variables $Y_t$, $t=1,\ldots, \nu$ take values, and to compute
the corresponding probabilities $\{f_t(j), \,j=1,\ldots, \ell_t\}$.
\par
The following algorithm implements the optimal policy.

\paragraph{Algorithm~1.}
\begin{itemize}
 \item[1.] Compute
 \[
  I_{t,k}(r)=\sum_{a=r}^{k-t+r} q(a) \frac{\binom{a-1}{r-1}\binom{k-a}{t-r}}{\binom{k}{t}},
\;\;\; r=1,\ldots, t;\;\;\;t=1,\ldots, k,
\]
where
\[
 k=\left\{ \begin{array}{ll}
             n, & {\rm Problem~(A1)},\\
             t,t+1,\cdots,N_{\max}(\hbox{ or } \tilde{N}_{\max}), & {\rm Problem~(A2)}.
            \end{array} \right.
\]

 We note that the computations
can be efficiently performed using the following recursive
formula:
for any reward function $q$
\begin{equation}\label{eq:mucci}
 I_{t,k}(r)=\frac{r}{t+1} I_{t+1,k}(r+1) + \Big(1-\frac{r}{t+1}\Big) I_{t+1,k}(r),\;\;\;r=1,\ldots,t;
\end{equation}
see \citeasnoun{GuZa} and \citeasnoun[Proposition~2.1]{mucci-a}.\\
Then compute
\begin{equation}\label{eq:U}
U_t(r)=\left\{ \begin{array}{ll}
             I_{t,n}(r), & {\rm Problem~(A1)},\\
             \sum_{k=t}^{\nu}\gamma_k I_{t,k}(r),  & {\rm Problem~(A2)}.
            \end{array} \right.
\end{equation}
\item[2.]
Find the distinct values $(y_t(1),\ldots, y_t(\ell_t))$ of the vector
$(U_t(1),\ldots,U_t(t))$, $t=1,\ldots, \nu$; here
$\ell_t$ is a number of the distinct points.
\item[3.] Compute
\[
 f_t(j)=\frac{1}{t}\sum_{r=1}^t
{\bf 1}\big\{ U_t(r)=y_{t}(j)\big\},\;\;\;\;j=1,\ldots,\ell_t; \;\;\;t=1,\ldots, \nu.
\]
\item[4.] Let $b_1=-\infty$, $b_2=\sum_{j=1}^{\ell_{\nu}} y_{\nu}(j) f_{\nu}(j)$.
\par
For $t=2, \ldots, \nu$ compute
\begin{equation}\label{eq:alg-recursion}
 b_{t+1}=\sum_{j=1}^{\ell_{\nu-t+1}} \big[\,b_t \vee y_{\nu-t+1}(j)\,\big]\,f_{\nu-t+1}(j).
\end{equation}
\item[5.] Output  $b_{\nu+1}$ and  $\tau_*=\min\{t\in\{1,\ldots, \nu\}: U_t(R_t)> b_{\nu-t+1}\}$.
In problems with random horizon, $\tau_*$ is the optimal stopping rule provided that
stopping occurred prior to termination of the observation process
due to horizon randomness.
\end{itemize}

\section{Solution of the sequential  selection problems}\label{sec:examples}

In this section we
revisit  problems (P1)--(P12)
discussed earlier  from the viewpoint of
the proposed framework. We refer to Section~\ref{sec:selection-problems} for
detailed description of these problems and related
literature.

\subsection{Problems with fixed horizon}
First we consider problems (P1)-(P5) with fixed horizon; in all these problems
$\nu=n$.
\subsubsection{Classical secretary problem}
For description of this problem and related references see Problem~(P1) in Section~\ref{sec:selection-problems}.
Here
$q(a)={\bf 1}\{a=1\}$, and
\[
U_t(r)=I_{t,n}(r)=\frac{t}{n}{\bf 1}\{r=1\},\;\;r=1,\ldots, t;\;\;
\;\;\ell_t=2, \;\;t=1,\ldots, n.
\]
The
random variable
$Y_t=(t/n) {\bf 1}\{R_t=1\}=\rP(A_{t,n}=1|R_t)$ takes two
different values
$y_t(1)=t/n$, $y_t(2)=0$
with probabilities
$f_t(1)=1/t$ and $f_t(2)=1- (1/t)$.
Then Step~4 of the Algorithm~1 takes the form:  $b_1=-\infty$, $b_2=1/n$,
\begin{eqnarray*}
\,\,\, b_{t+1}
=
b_t+ \Big(\frac{1}{n}-\frac{b_t}{n-t+1}\Big) {\bf 1}\Big\{b_t<\frac{n-t+1}{n}\Big\},\,\,t=2,\dots,n.
\end{eqnarray*}
The optimal policy is to stop the first time instance $t$ such that $Y_t> b_{n-t+1}$, i.e.,
\[
 \tau_*=\min\Big\{1\leq t\leq n: \frac{t}{n}\,{\bf 1}\{R_t=1\}> b_{n-t+1}\Big\},
\]
which coincides with  well known results. 

\subsubsection{Selecting one of $k$ best alternatives}\label{sec:Gusein-Zade}
This setting is stated as Problem~(P2) in Section~\ref{sec:selection-problems}.
In this problem $q(a)={\bf 1}\{a\leq k\}$ with some $k\leq n$.
We will assume here that $k\geq 2$; the case $k=1$ was treated above.
\par
We have
\begin{equation}\label{eq:U-U}
 U_t(r)=\left\{
\begin{array}{ll}
0, & k+1\leq r \leq t,
\\*[4mm]
\sum_{a=r}^{(n-t+r)\wedge k} \frac{\binom{a-1}{r-1}\binom{n-a}{t-r}}{\binom{n}{t}},
&  1\leq r \leq k,
\end{array}
\right.\;\;\;\;\;t=1,\ldots, n.
\end{equation}
It is easily checked that for $q(a)={\bf 1}\{a\leq k\}$ one has
\begin{equation}\label{eq:initial-condition}
 U_n(r)=\left\{\begin{array}{ll}
                1, & r=1,\ldots,k\\
0, & r=k+1,\ldots,n.
               \end{array}
\right.
\end{equation}
Using this formula together with the recursive relationship (\ref{eq:mucci})
we can determine the structure of  vector $U_t:=(U_t(1),\ldots, U_t(t))$ for each
$t=1,\ldots, n$,
and compute $\{y_t(j)\}$ and $\{f_t(j)\}$. Specifically, the following facts
are  easily verified.
\begin{itemize}
\item[(a)] Let $n-k+2\leq t\leq n$. Here
vector
$U_{t}$ has the following structure:
the first 
$t+k-n$ components are ones, the next
$
n-t$ components are distinct numbers
in $(0,1)$ which are given in
(\ref{eq:U-U}), and the last $
t-k$ components are zeros. Formally,
if $n-k+2\leq t \leq n-1$ and $k>2$ then
we have
\[
 U_t(j)=\left\{\begin{array}{ll}
                1, & j=1,\ldots,k-n+t,\\
                \in (0,1), & j=k-n+t+1,\ldots, k,\\
                0, & j=k+1,\ldots, t,
               \end{array}
\right.
\]
Note that if $k=2$ the regime reduces to $t=n$; therefore
if $k=2$ or $t=n$
then $U_n$ is given by (\ref{eq:initial-condition}).
These facts imply the following expressions
for $\{y_t(j)\}$ and $\{f_t(j)\}$:
\begin{eqnarray}\label{eq:gusein-(c1)}
 \ell_t=n-t+2;\;\;\;
y_t(j)=\left\{\begin{array}{ll}
               1, & j=1,\\
   U_t(k-n+t+j), & j=2,\ldots,n-t+1,\\
0, & j=n-t+2,
              \end{array}
\right.
\end{eqnarray}
and
\begin{eqnarray}\label{eq:gusein-(c2)}
f_t(j)=\left\{\begin{array}{ll}
               1- (n-k)/t, & j=1,\\
1/t, & j=2,\ldots, n-t+1,\\
1- k/t, & j=n-t+2.
              \end{array}
\right.
\end{eqnarray}
If $t=n$ then
\[
\ell_t=2, \;y_n(1)=1,\; y_n(2)=0,\; f_n(1)=k/n,\; f_n(2)=1-k/n.
\]
\item[(b)]
If $k+1\leq t\leq n-k+1$ then the set $\{U_t(1),\ldots,U_t(t)\}$ contains $k+1$ distinct values:
$U_t(1), \ldots, U_t(k)$ are positive distinct, and $U_t(k+1)=\cdots=U_t(t)=0$.
Therefore
\begin{equation}\label{eq:gusein-(b)}
 \ell_t=k+1;\;\;\;y_t(j)=\left\{\begin{array}{ll}
                                 U_t(j), & j=1,\ldots,k\\
0, & j=k+1;
                                \end{array}
\right.\;\;\;f_t(j)=\left\{\begin{array}{ll}
                            1/t, & j=1,\ldots,k,\\
1- k/t, & j=k+1.
                           \end{array}
\right.
\end{equation}
\item[(c)] If $1\leq t\leq k$ then all the values $U_t(1),\ldots, U_t(t)$ are positive and distinct.
Thus
\begin{equation}\label{eq:gusein-(a)}
 \ell_t=t;\;\;\;y_t(j)=U_t(j),\;j=1,\ldots,t;\;\;\;f_t(j)=\frac{1}{t},\;\;j=1,\ldots,t.
\end{equation}
\end{itemize}
\par
In our implementation we compute $U_t(j)$ for $t=1,\ldots, n$ and $j=1,\ldots, t$ using (\ref{eq:initial-condition}) and (\ref{eq:mucci}). Then
$\{y_t(j)\}$, $\{f_t(j)\}$ and the sequence $\{b_t\}$ are easily calculated from
(\ref{eq:gusein-(c1)})--(\ref{eq:gusein-(a)}) and
(\ref{eq:alg-recursion}) respectively.
\par
Table~\ref{tab:gusein} presents exact values of
the optimal probability $P(n,k)=b_{n+1}$ and the
expected time until stopping $E(n,k)=\rE (\tau_*)$ normalized by $n$
for different values of
$k$ and $n$. We are not aware of works that  report exact results  for
general $k$ and $n$ as presented in Table~\ref{tab:gusein}.
These results should be compared to
the asymptotic values of $1-P(n,k)$ as $n\to\infty$
computed
in \citeasnoun[Table~1]{FrSa}
for a range of values of $k$. The comparison shows that the approximate
values in \citeasnoun{FrSa} are in a good agreement with the
exact values of Table~\ref{tab:gusein}.
For instance, for $n=100$ the approximate values coincide with the exact ones up to the third digit
after the decimal point.
\par
It is worth noting that
the optimal policy developed by
\citeasnoun{GuZa}
is expressed in terms of
of relative ranks.
In contrast, our policy is
expressed via the random variables $Y_t=U_t(R_t)$, and it is
memoryless threshold in terms of $\{Y_t\}$.
This  allows to efficiently compute the
distribution of the optimal stopping time, and, in particular, the
expected time until stopping.
The value of $E(n,k)$ is computed using formula
(\ref{ex:ExpTimeToStop}) combined with (\ref{eq:distr-Y}) and (\ref{eq:U-U})--(\ref{eq:gusein-(a)}).
The presented numbers agree with asymptotic results of \citeasnoun{yeo} proved for $k=2,3$ and $5$.
\begin{table}
\begin{center}
{\tiny
 \begin{tabular}{||c|c|c|c||c|c|c|c||c|c|c|c||}\hline\hline
  $n$ & $k$ & $P(n,k)$ &$E(n,k)/n$ & $n$ & $k$ & $P(n,k)$&$E(n,k)/n$& $n$ & $k$ & $P(n,k)$&$E(n,k)/n$\\
\hline\hline
     100 & 2 &  0.57956 &0.68645&500 &  2   & 0.57477&0.68886& 1,000 & 2 &  0.57417&0.68966\\
         & 5 & 0.86917  &0.60871 &     &  5  &   0.86211&0.60921&     &  5 &  0.86123&0.60988  \\
         & 10& 0.98140  &0.54236  &    &  10  & 0.97754&0.54454 &    &  10 & 0.97703&0.54434  \\
         & 15 &  0.99755 &0.50428 &     &  15 &  0.99627&0.50845  &     &  15&  0.99609&0.50893\\
         \hline\hline
   5,000 & 2  &  0.57369 &0.68931 &10,000 & 2 & 0.57363&0.68927 &50,000 &2 & 0.57358&0.68923  \\
         &  5&  0.86052& 0.61015 &         &5&  0.86043&0.61014 &       & 5 &0.86036&0.61018 \\
         & 10 &0.97663 & 0.54499 &        &10& 0.97658&0.54496  &      & 10&0.97654& 0.54500 \\
         & 15 & 0.99594 & 0.50943&        &15&  0.99592&0.50947 &      &15 &0.99591& 0.50950\\
         \hline\hline
 \end{tabular}
}
\caption{Optimal probabilities $P(n,k)$ and the normalized expected time elapsed until stopping $E(n,k)/n$ for selecting one of the $k$ best values.
\label{tab:gusein}}
\end{center}
\end{table}

\subsubsection{Selecting the $k$-th best alternative}
This setting is discussed in Section~\ref{sec:selection-problems} as problem~(P3).
In this problem  $q(a)={\bf 1}\{a=k\}$, $k\geq 2$.
Similarly to the Gusein--Zade stopping problem, here we have three different regimes
that define explicit relations for $\{U_t(r)\}$, $\{y_t(j)\}$ and $\{f_t(j)\}$.
\begin{enumerate}
\item[(a)] Let
$1\leq t \leq k$; then
\begin{align*}
&
{\displaystyle
U_t(r)=
\frac{\tbinom{k-1}{r-1}\tbinom{n-k}{t-r}}{\tbinom{n}{t}},\,\,r=1,\ldots,t.
}
\end{align*}
All values of ${\displaystyle U_t(1),\ldots,U_{t}(t)}$ are positive and distinct. Thus
\begin{equation}\label{eq:postdoc(a)}
 \ell_t=t,\,\,\; y_t(j)=U_{t}(j),\,\,\;f_{t}(j)=\frac{1}{t},\,\, 1\leq j \leq t.
\end{equation}
\item[(b)]
If $k+1\leq t \leq n-k+1$ then
\begin{align*}
&U_t(r)=
\left\{
\begin{array}{ll}
\frac{\tbinom{k-1}{r-1}\tbinom{n-k}{t-r}}{\tbinom{n}{t}}, & 1\leq r\leq k,
\\
0, & k+1\leq r \leq t.
\end{array}
\right.
\end{align*}
The set $\{U_t(1),\ldots,U_t(t)\}$ contains $k+1$ distinct values:
$U_t(1), \ldots, U_t(k)$ are positive distinct, and $U_t(k+1)=\cdots=U_t(t)=0$.
Therefore,
\begin{equation}\label{eq:postdoc(b)}
 \ell_t=k+1;\;\;\;y_t(j)=\left\{\begin{array}{ll}
                                 U_t(j), & j=1,\ldots,k\\
0, & j=k+1;
                                \end{array}
\right.\;\;\;f_t(j)=\left\{\begin{array}{ll}
                            1/t, & j=1,\ldots,k,\\
1- k/t, & j=k+1.
                           \end{array}
\right.
\end{equation}
\item[(c)] Let
$n-k+2\leq t \leq n$; then
the sequence $\{U_t(r)\}$ takes the following values
\[
 U_t(r)=\left\{\begin{array}{ll}
                0, & r=1,\ldots,t-n+k-1,\\
                \frac{\tbinom{k-1}{r-1}\tbinom{n-k}{t-r}}{\tbinom{n}{t}}, & r=t-n+k,\ldots, k,\\
                0, & r=k+1,\ldots, t.
               \end{array}
\right.
\]
Therefore,
\begin{equation}\label{eq:postdoc(c1)}
 \ell_t=n-t+2;\;\;\; y_t(j)=\left\{\begin{array}{ll}
                                0, &j=1\\
                                 U_t(t-(n-k)-2+j), & j=2,\ldots, n-t+2,\\
                               \end{array} \right.
\end{equation}
and, correspondingly,
\begin{equation}\label{eq:postdoc(c2)}
 f_t(j)= \left\{\begin{array}{ll}
                 (2t-n-1)/t, & j=1,\\
                 1/t, & j=2,\ldots, n-t+2.
                \end{array}
\right.
\end{equation}
\end{enumerate}
\par
Table~\ref{tab:postdoc} presents optimal probabilities of selecting $k$th best alternative
for a range of  $k$ and $n$.  In the specific case of $k=2$ \citeasnoun{Rose}
showed that the optimal stopping rule is
\[
\tau_*=\min\Big\{ \{t\geq \lceil n/2\rceil: R_t=2\}\cup \{n\} \Big\},
\]
and the optimal probability is $P(n,2)=\frac{n+1}{4n}$ if $n$ is odd. The results for $k=2$
in Table~\ref{tab:postdoc} are in full agreement with this formula. The table also presents
numerical computation of optimal values in the problem of selecting the median value; see
\citeasnoun{Rose-2} who proved that $\lim_{n\to\infty} V^*_n(q_{\rm pd}^{((n+1)/2)})=0$.
\begin{table}
\begin{center}
{\tiny
 \begin{tabular}{||c|c|c||c|c|c|c|c||c|c|c|c||}\hline\hline
  $n$ & $k$ & $P(n,k)$&$E(n,k)/n$ & $n$ & $k$ & $P(n,k)$&$E(n,k)/n$& $n$ & $k$ & $P(n,k)$&$E(n,k)/n$\\
\hline\hline
     101 & 2 & 0.25247&0.82995  &501 &  2   &  0.25050 &0.75466& 1,001 & 2 & 0.25025&0.74984\\
         & 5 & 0.19602&0.78968  &    &  5  &  0.19281&0.78890&     &  5 &  0.19241&0.78896  \\
         & 10& 0.15962&0.84827   &   &  10  & 0.15506&0.84508&      &  10 &0.15451&0.84517  \\
         & 50 & 0.11467&0.86699   &   &  250 & 0.06876&0.91156 &      & 500&  0.05504&0.92688\\
         \hline\hline
   5,001 & 2  &0.25005 & 0.84527& 10,001 & 2 &0.25002 &0.75453 & 50,001 &2 &  0.25000&0.83830\\
         &  5& 0.19210 &0.78896&           &5& 0.19206&0.78891 &       & 5 &0.19203&0.78891\\
         & 10 &0.15450 &0.84478  &        &10& 0.15402&   0.84477&     & 10&0.15397&0.84477 \\
         & 2,500 & 0.03265&0.95443 &         &5,000& 0.02603&0.96320&       &25,000 &0.01533&0.97787 \\
         \hline\hline
 \end{tabular}
}
\caption{Optimal probabilities $P(n,k)$ and the normalized expected time elapsed until
stopping $E(n,k)/n$ for  selecting the $k$-th best
alternative  computed using~(\ref{eq:postdoc(a)})--(\ref{eq:postdoc(c2)}).
\label{tab:postdoc}}
\end{center}
\end{table}

\subsubsection{Expected rank type problems}
In this section we consider problems~(P4) and (P5) discussed in Section~\ref{sec:selection-problems}.

\paragraph{Expected rank minimization.}
Following (\ref{eq:chow-rob}) we
consider the problem of minimization of $\rE q(A_{\tau,n})$, where
$q(a)=-a$.
It is well known that $\rE\big[A_{t,n}|R_t=r\big] = (n+1)r/(t+1)$; therefore
for $t=1,\ldots, n$
\[
 U_t(r)=I_{t,n}(r)= \rE [q(A_{t,n})| R_t=r]= -\rE
 [A_{t,n}|R_t=r]=- \frac{(n+1)r}{t+1}, \;\;r=1,\ldots, t.
\]
In this setting
\[
\ell_t=t,\;\;\forall t;\;\;\; y_t(j)=U_t(j)=-\frac{n+1}{t+1}j,\;\;j=1,\ldots t;\;\;\;
f_t(j)=\frac{1}{t},\;\;\;\ \forall j=1,\ldots,t.
\]
Substitution to (\ref{eq:rec-b}) yields $b_1=-\infty$, $b_2=-\frac{1}{2}(n+1)$,
\begin{equation}
\label{eq:Chow}
 b_{t+1}= \frac{1}{n-t+1} \sum_{j=1}^{n-t+1}
 \Big[b_t \vee \Big(-\frac{n+1}{n-t+2}j\Big)\Big],\;\;\;
t=2,\ldots,n.
\end{equation}
Straightforward calculation shows that  (\ref{eq:Chow}) takes form
\[
 b_{t+1}= b_t - \frac{1}{n-t+1}  \bigg[\frac{n+1}{n-t+2}\frac{j_t(j_t+1)}{2}
 +j_t b_t \bigg],\;\;\;
t=2,\ldots,n.
\]
where
$j_t:=\lfloor -b_t \frac{n-t+2}{n+1}\rfloor$.
The optimal policy is to stop the first time instance
$t$ such that $Y_t > b_{n-t+1}$, i.e.,
\begin{align*}
 \tau_{*}=\min\Big\{1\leq t\leq n: -\frac{n+1}{t+1}R_t> b_{n-t+1}\Big\}
 =\min\Big\{1\leq t\leq n: R_t\, \leq \, j_{n-t+1}  \Big\}.
\end{align*}
Then according to (\ref{eq:chow-rob}) the optimal value of the problem equals to
$-b_{n+1}$.
We note that the derived recursive procedure coincides with the one of \citeasnoun{chow}, and
the calculation for $n=10^6$ yields the optimal value $3.86945\ldots$

\paragraph{Expected squared rank minimization.}
This problem was posed in \citeasnoun{Robbins-91}, and to the best of
our knowledge, it was not solved to date. We show that the proposed unified
framework can be used in order to compute efficiently the optimal policy and its value.
\par
In this setting  $U_{t}(r)=I_{t,n}(r)$, and
the reward is given by $q(a)=-a^2$.
It is well known that
\[
\rE\big[A_{t,n}(A_{t,n}+1)\cdots (A_{t,n}+k-1)\,|\,R_t=r\big]  = \frac{(n+1)\cdots (n+k)}{(t+1)\cdots (t+k)}
 r\cdots (r+k-1);
\]
see, e.g., \citeasnoun{Robbins-91}.
Therefore we put
\begin{eqnarray*}
 U_t(r) = -\rE (A_{t,n}^2|R_t=r)
 =-
 \frac{(n+1)(n+2)}{(t+1)(t+2)}r  \Big( r +
 \frac{n-t}{n+2}\Big).
\end{eqnarray*}
In this case
\begin{eqnarray*}
&& \ell_t=t, \;\;y_t(j)=U_t(j)=-\frac{(n+1)(n+2)}{(t+1)(t+2)}j   \Big(j +
 \frac{n-t}{n+2}\Big),\;\;f_t(j)=\frac{1}{t},\;\;\;j=1,\ldots, t.
\end{eqnarray*}
Substituting this to (\ref{eq:rec-b}) we obtain the following  recursive relationship: $b_1=-\infty$,
$b_2=-\tfrac{1}{6}(n+1)(2n+1)$,
\begin{eqnarray*}
b_{t+1}&=&\frac{1}{n-t+1}\sum_{j=1}^{n-t+1}
\bigg\{b_t\vee \Big[-\frac{(n+1)(n+2)}{(n-t+2)(n-t+3)}j\Big(j+\frac{t-1}{n+2}\Big)\Big]\bigg\}.
\label{eq:bt}
\end{eqnarray*}
Denote
$j_t:=\max\{1\leq j\leq n-t+1: b_t \leq -j^2C_{n,t}-jD_{n,t}\}$,
where
\[
C_{n,t}=\frac{(n+1)(n+2)}{(n-t+2)(n-t+3)},\,\,\,\,D_{n,t}=\frac{(t-1)(n+1)}{(n-t+2)(n-t+3)}.
\]
Then
\begin{eqnarray*}
j_t=\max\bigg\{1\leq j\leq n-t+1: j \leq \frac{1}{2C_{n,t}}
\Big(-D_{n,t}+
\sqrt{D^2_{n,t}-4C_{n,t}b_t}\,\Big)  \bigg\}
\\*[2mm]
=
\Big\lfloor \frac{1}{2C_{n,t}}\Big(-D_{n,t}+
\sqrt{D^2_{n,t}-4C_{n,t}b_t}\,\Big) \Big\rfloor.
\end{eqnarray*}
With this notation we have
$b_1=-\infty$, $b_2=-\frac{1}{6}(n+1)(2n+1)$, and for $t=2,\ldots, n$
\begin{align}\label{eq:A2}
b_{t+1}= \frac{1}{n-t+1}  \bigg[-\frac{1}{6}j_t(j_t+1)(2j_t+1)C_{n,t}
-\frac{1}{2}j_t(j_t+1)D_{n,t} +(n-t+1-j_t) b_t \bigg].
\end{align}
The optimal policy is to stop the first time instance $t$ such that $Y_t> b_{n-t+1}$
which is equivalent to
\begin{align*}
&
 \tau_{*}=\min\Big\{1\leq t\leq n: R_t\leq j_{n-t+1}  \Big\}.
\end{align*}
Table~\ref{tab:squared} presents optimal values $V_*(n):=\rE A_{\tau_*,n}^2$
computed with recursive relation (\ref{eq:A2})
for different $n$.
\begin{table}
\begin{center}
{\small
 \begin{tabular}{||c||c|c|c|c|c|c|}\hline\hline
  $n$ & 100 & 250 & 500 & 750 & 1,000 & 2,500\\
\hline
$V_*(n)$ &23.70663   &26.49268  &27.66697  &28.10937 & 28.34466 & 28.80553   \\
\hline\hline
$n$  & 5,000 & 10,000 & 20,000 & $10^5$ & $10^6$ & $10^8$
\\
\hline
 $V_*(n)$  &28.97697 & 29.06969 &29.11944  & 29.16302 & 29.17431 & 29.17579\\
 \hline\hline
 \end{tabular}
}
\caption{Optimal values of $V_*(n):=\rE A_{\tau_*,n}^2$ computed using (\ref{eq:A2}).
\label{tab:squared}}
\end{center}
\end{table}
\subsection{Problems with random horizon}
This section demonstrates how to apply the proposed framework for solution of
selection problems with a random horizon.
In these problems we apply Algorithm~1 with $\nu$
being the maximal horizon length $N_{\max}$, provided that $N_{\max}$ is finite,
or with sufficiently large horizon $\tilde{N}_{\max}$ if $N_{\max}$ is infinite. Moreover,
$U_t(r)=J_t(r)$,
where $\{J_t(r)\}$ is given by (\ref{eq:tilde-u}).
\par
Recall that in all problems with random horizon
the selection may not be made by the time 
the observation process terminates.
However, Theorems~\ref{th:random} and~\ref{th:random-reduction} show that
as long as the observation process proceeds, the optimal stopping rule is identical to
the one in the setting with fixed horizon $N_{\max}$
and random variables $Y_t:=U_t(R_t)$, $t=1, \ldots, N_{\max}$, where  $U_t(\cdot)$
is defined in (\ref{eq:U}).
In the subsequent discussion of specific problem instances with random horizon
we use this fact without further mention.
\subsubsection{Classical secretary problem with random horizon}
This is Problem~(P5) of Section~\ref{sec:selection-problems} where
$q(a)={\bf 1}\{a=1\}$; therefore
\begin{eqnarray*}
 && I_{t,k}(r)=\rP(A_{t, k}=1\,|\,R_t=r)= \frac{t}{k} {\bf 1}\{r=1\},\;\;k\geq t,
 \\
 && U_t(r)=J_t(r)= \sum_{k=t}^{N_{\max}} \gamma_k I_{t,k}(r)=
 t\, {\bf 1}\{r=1\} \sum_{k=t}^{N_{\max}} \frac{\gamma_k}{k}~.
\end{eqnarray*}
Note that if $N_{\max}=\infty$ then condition (\ref{eq:u-finite}) is trivially fulfilled
since
\[
t \sum_{k=t}^{\infty} \frac{\gamma_k}{k} \leq \sum_{k=t}^\infty \gamma_k \leq 1.
\]
\par
The
random variables
$Y_t=U_t(R_t)={\bf 1}(R_t=1)\,t\,\sum_{k=t}^{\nu}\gamma_k/k$ take two
different values $y_t(1)=t\sum_{k=t}^{\nu}\gamma_k/k$ and
$y_t(2)=0$ with corresponding probabilities
$f_t(1)=1/t$ and $f_t(2)=1- 1/t$.
Substituting these values in
(\ref{eq:alg-recursion}) we obtain $b_1=-\infty$, $b_2=\gamma_{\nu}/\nu$, and
for $t=2,\ldots, \nu$
\begin{eqnarray}
b_{t+1}= b_t+\bigg(\sum_{k=\nu-t+1}^{\nu}\frac{\gamma_k}{k}-\frac{b_t}
{\nu-t+1}
\bigg){\bf 1}\bigg\{b_t< (\nu-t+1)
\sum_{k=\nu-t+1}^{\nu}\frac{\gamma_k}{k}\bigg\}.
\label{eq:class-RRusp}
\end{eqnarray}
The optimal policy is to stop at time $t$ if $Y_t > b_{\nu-t+1}$, i.e.,
\begin{align}
\label{eq:tau-csp-random}
\tau_{\ast}=
\min\;\Big\{t=1,\ldots,\nu:
\,{\bf 1}\{R_t=1\}\,t\,\sum_{k=t}^{\nu}\frac{\gamma_k}{k}> b_{\nu-t+1}
\Big\}.
\end{align}
\par
\citeasnoun{PS1972}
investigated  the structure of  optimal stopping rules  and
showed that, depending on the distribution of $N$,  the stopping region can involve  several ``islands,'' i.e., it  can be a union of disjoint subsets
of $\{1,\ldots, N_{\max}\}$.
Note that (\ref{eq:tau-csp-random}) determines the stopping region automatically.
Indeed, it is optimal to stop only at those $t$'s that satisfy
$t \sum_{k=t}^\nu \gamma_k/k > b_{\nu-t+1}$.
We apply the stopping rule (\ref{eq:class-RRusp})--(\ref{eq:tau-csp-random})
for two examples of distributions of $N$. In the first example $N$ is assumed to be
uniformly distributed on the set $\{1,\ldots, N_{\max}\}$. As it is known, in this case
the optimal stopping region has only one ``island.'' The second example illustrates a setting in which
the stopping region has more than one ``island.''
\par\medskip
1. {\em Uniform distribution.}
In this case $\nu=N_{\max}$, $\gamma_k=1/N_{\max}$, $k=1,\ldots, N_{\max}$.
It was shown in
 \citeasnoun{PS1972} that
 the optimal stopping region in this problem has one ``island,'' i.e.,
 the optimal policy selects
 the first best member appearing in the range $\{k_n, \ldots, n\}$.
The
recursive relation (\ref{eq:class-RRusp}) with $\gamma_k=1/N_{\max}$, $k=1,\ldots, N_{\max}$
yields the  optimal values $V_*(N_{\max}):= \rP\{A_{\tau_*, N}=1, \tau_*\leq N\}$
given in Table~\ref{tab:uniform-d}.
The second line of Table~\ref{tab:uniform-d} presents the normalized expected time
until stopping $E_*(N_{\max}):= \rE(\tau_*\wedge N_{\max})/ N_{\max}$ computed using (\ref{eq:tau-min}), (\ref{eq:111}) and (\ref{eq:222}).
For comparison, we also give the normalized expected time elapsed until stopping
$E_*(n):=\rE \tau_*/n$ for the optimal stopping rule in the classical secretary problem (see the third line of the table). These numbers are calculated using (\ref{ex:ExpTimeToStop}). As expected, $E_*(N_{\max})$
is significantly smaller than $E_*(n)$;  the optimal rule is more cautious when the horizon is random.

\begin{figure}
\[
    \begin{array}{cc}
    \hspace{-7mm} \includegraphics[scale=0.21]{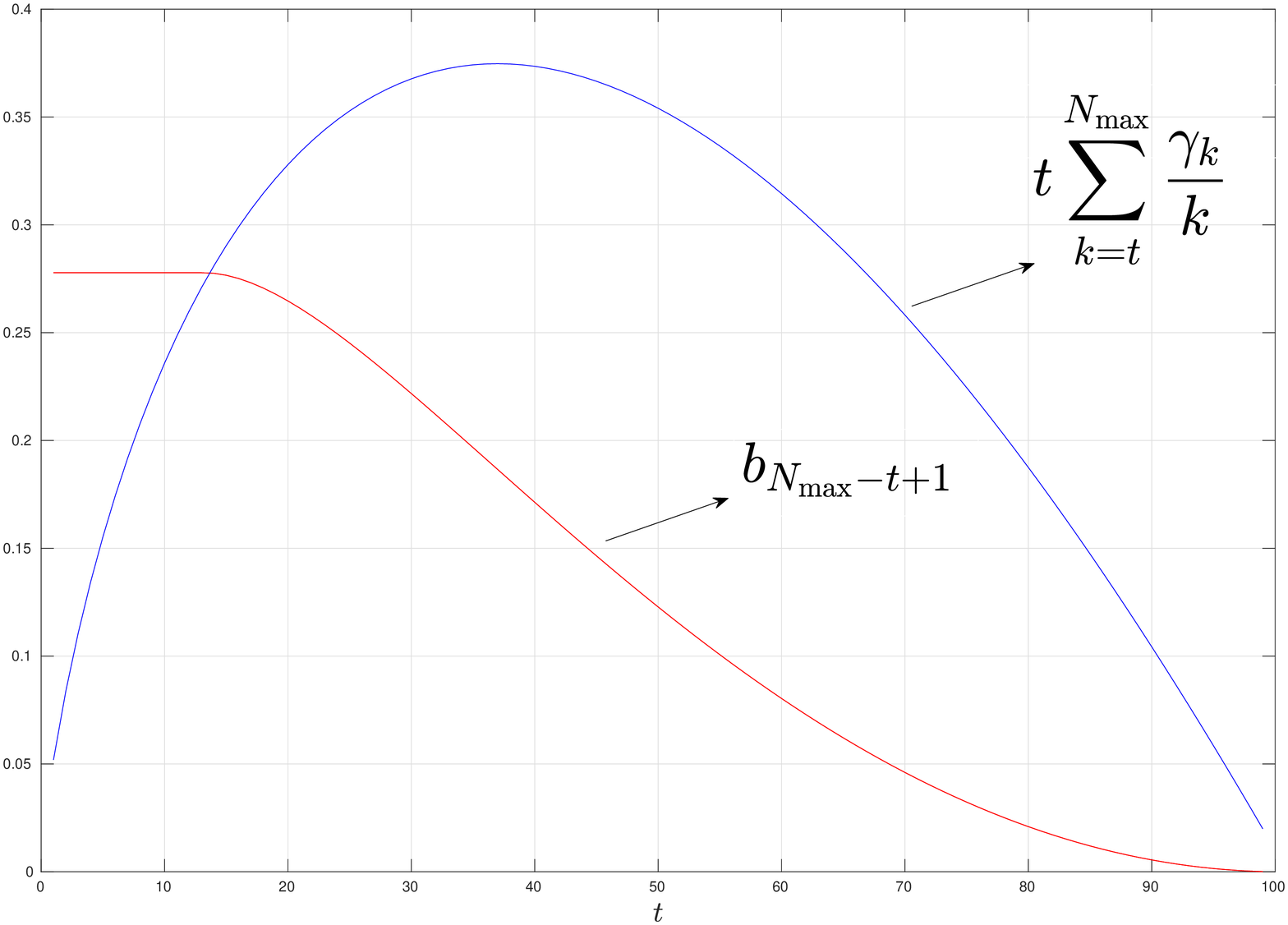} & \hspace{-5mm}
    \includegraphics[scale=0.235]{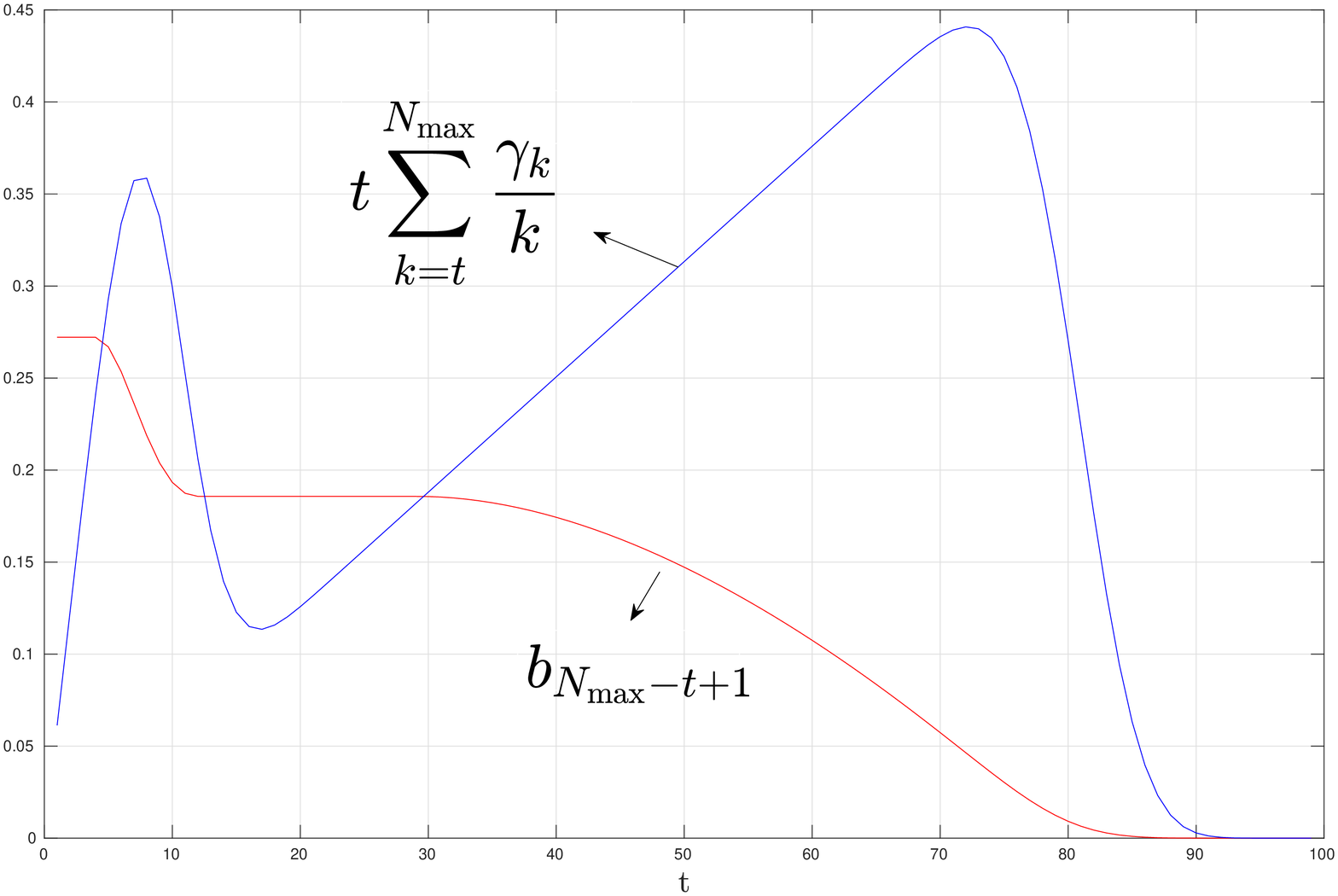}
    \\
    {\rm (a)} & {\rm (b)}
    \end{array}
\]
    \caption{The graphs of sequences $\{b_{N_{\max}-t+1}\}$ and $\{t\sum_{k=t}^{N_{\max}} \gamma_k/k\}$
    for different distributions of $N$: (a) the uniform distribution; (b)  the mixture of two zero--inflated binomial distributions.
    \label{fig:UnRank}}
\end{figure}

\begin{table}
\begin{center}
{\small
\begin{tabular}{||c|c|c|c|c|c|c|c|c||} \hline\hline
$N_{\max} \,|\, n $&10&20&40&60&80&$10^2$&$10^3$&$10^5$\\ \hline
$V_*(N_{\max})$&0.35145& 0.30760& 0.28889& 0.28260& 0.27949&0.27779&0.27137&0.27068\\
$E_*(N_{\max})$&0.29290& 0.26227   &0.280651 & 0.28605  & 0.27410   & 0.27410  & 0.27995 & 0.27983\\
$E_*(n)$& 0.61701 &  0.73421  & 0.75074   & 0.73988   & 0.73436   & 0.74104 &0.73620  &0.73576\\
\hline\hline
\end{tabular}
}
\caption{Optimal values $V_*(N_{\max}):=\rP\{A_{\tau_*,N}=1, \tau_*\leq N\}$ for a uniformly  distributed
horizon length $N$, normalized expected times until stopping $E_*(N_{\max})$ and
$E_*(n)$ for random and fixed horizons.
\label{tab:uniform-d}}
\end{center}
\end{table}
%
%
It was also shown in \citeasnoun{PS1972} that
$\lim_{N_{\max}\to \infty}V_*(N_{\max})=2e^{-2}=0.27067\ldots$. Note that the  numbers in
Table~\ref{tab:uniform-d} are in full agreement with these results.
Figure~\ref{fig:UnRank}(a) displays the sequences $\{b_{N_{\max}-t+1}\}$ and
$\big\{t\sum_{k=t}^{N_{\max}} \gamma_k/k\big\}$ for the uniform distribution for $N_{\max}=100$.
Note the  stopping region  is the set of $t$'s where the blue curve is above the red curve. Thus,  there is
only one ``island'' in this case.
\par\medskip

\par\medskip
2. {\em Mixture of two zero--inflated binomial distributions.} Here we assume that
the distribution $G_N$ of $N$ is  the mixture: $G_N(x)= \frac{1}{2}H_1(x)+ \frac{1}{2} H_2(x)$, where
$H_i(x) = \rP(X_i\leq x| X_i\geq 1)$, $i=1,2$,  and  $X_1\sim {\rm Bin}(50, 0.2)$, $X_2\sim {\rm Bin}(100, 0.8)$.
In other words, for $k=1,\ldots, 100$
\[
 \gamma_k=\rP(N=k)=\frac{1}{2}
 \binom{50}{k} \Big(\frac{1}{4}\Big)^k \frac{(0.8)^{50}}{1-(0.8)^{50}}
 +\frac{1}{2} \binom{100}{k} 4^k \frac{(0.2)^{100}}{1-(0.2)^{100}}.
\]
The optimal stopping rule is given by (\ref{eq:class-RRusp})--(\ref{eq:tau-csp-random})
with $\{\gamma_k\}$ indicated above.
Figure~\ref{fig:UnRank}(b) displays
the graphs of the sequences $\{b_{N_{\max}-t+1}\}$ and
$\big\{t\sum_{k=t}^{N_{\max}} \gamma_k/k\big\}$. It is clearly seen that in this setting
the stopping region is a union of two disjoint sets of subsequent integer numbers. These sets
correspond to the indices where the graph of
$\big\{t\sum_{k=t}^{N_{\max}} \gamma_k/k\big\}$
is above the graph of $\{b_{N_{\max}-t+1}\}$.
The stopping region can be easily
identified from given formulas.
\subsubsection{Selecting one of $k$ best alternatives with random horizon}
This is Problem~(P6) of Section~\ref{sec:selection-problems}; here $q(a)={\bf 1}\{a\leq k\}$.
Algorithm~1 is implemented
similarly to Problem~(P2). First, values
$I_{t,k}(r), k=1,\ldots,N_{\max}, t=1,\ldots,k, r=1,\ldots,t$ are
calculated using the recursive formula \eqref{eq:mucci} along with the boundary condition
(\ref{eq:initial-condition}). Then, using \eqref{eq:U}, we compute
$U_t(1),\ldots,U_t(t)$ for $t=1,\ldots,N_{\max}$,
and find  the distinct values $y_t(1),\ldots,y_t(\ell_t)$ of the vector $(U_t(1),\ldots,U_t(t))$ for all
$t=1,\ldots,N_{\max}$.
Finally, the  sequence $\{b_{t}\}$ is found from \eqref{eq:alg-recursion}.
The optimal policy is to stop the first time instance $t$ such that $Y_t=U_t(R_t)>b_{n-t+1}$
provided that the observed relative rank is different from zero; otherwise, the selection process terminates by the
problem horizon $N$.
The optimal value of the problem is  $P(N_{\max},k):=\rP\{A_{\tau_{*}, N}\leq k, \tau_{*}\leq N\}=b_{N_{\max}+1}$.
We apply this algorithm for two different examples:
a uniform horizon distribution, and
a~U--shaped distribution.
The second example demonstrates that  the optimal stopping region can have
``islands''  in the terminology of \citeasnoun{PS1972}.
\par\medskip
1. {\em Uniform distribution.}
In this case $\gamma_k=1/N_{\max}$, $k=1,\ldots, N_{\max}$.
Table \ref{tab:guseinUniform} presents exact values of the optimal probability $P(N_{\max},k)$. For $k=1$ the values
of $P(N_{\max},1)$ are in agreement with the values of Table \ref{tab:uniform-d} and also with the asymptotic
value obtained by \citeasnoun{PS1972}, $\lim_{N_{\max}\to \infty}P(N_{\max},1)=2e^{-2}=0.27067\ldots$.
For $k=2$ the values of $P(N_{\max},2)$ are in the agreement with the values of Table 1 in \citeasnoun{KT2003} and also with the asymptotic value obtained there, $\lim_{N_{\max}\rightarrow \infty }P(N_{\max},2)\approx 0.4038$.

\begin{table}
\begin{center}
{\small
 \begin{tabular}{||c|c|c|c||c|c|c|c||c|c|c|c||}\hline\hline
  $N_{\max}$ & $k$ & $P(N_{\max},k)$ & $N_{\max}$ & $k$ & $P(N_{\max},k)$& $N_{\max}$ & $k$ & $P(N_{\max},k)$\\
\hline\hline
     100 & 1 & 0.27779  &500&    1   &0.27208 & 1,000 & 1& 0.27137\\
         & 2 & 0.41506  & &     2    &0.40606 &       & 2 &0.40494 \\
         & 5 & 0.61788  & &     5    &0.60351  &      & 5 &0.60174  \\
         & 10& 0.75150  &  &    10   &0.73303 &       & 10&0.73078 \\
         & 15 &0.81474  &      &15   &0.79415 &       & 15&0.79161  \\
         \hline\hline
   5,000 & 1  &0.27081    &10,000 & 1   & 0.27074   &50,000 &1  &0.27068  \\
         &  2 &0.40405 &          &2    & 0.40394         && 2 & 0.40385 \\
         & 5  &0.60033 &          &5    & 0.60015         && 5 &0.60001  \\
         & 10 &0.72899  &         &10   & 0.72877         &&10 &0.72859 \\
         & 15 &0.78961  &         &15   & 0.78936         &&15 &0.78916 \\
         \hline\hline
 \end{tabular}
}
\caption{Optimal values $P(N_{\max},k):=P(A_{\tau_{*}}\leq k, \tau_{*}\leq N)$ for a uniformly distributed horizon length $N$.}
\label{tab:guseinUniform}
\end{center}
\end{table}
\par
\begin{figure}[h]
\begin{center}
\includegraphics[scale=0.40]{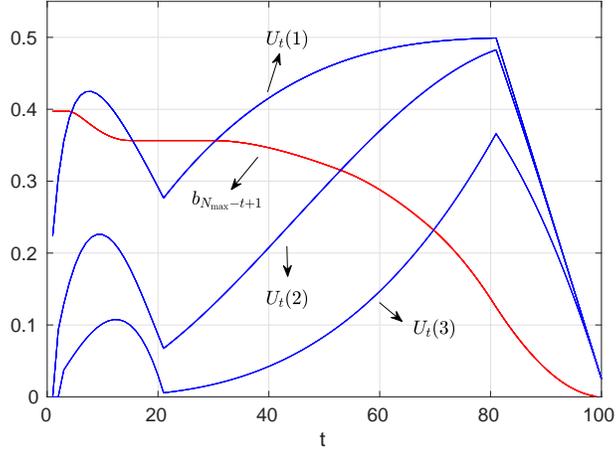}
\caption{The graphs of sequences $\{b_{N_{\max}-t+1}\}$, $\{U_t(r)\}$, $t=1,2,3$
    for the U-shape distribution distribution of $N$ defined in (\ref{eq:U-dist}).
    \label{fig:Ushape}}
\end{center}
\end{figure}

\medskip
2. {\em U-shaped distribution.}
In this example we let  $N_{\max}=100$,
\begin{equation}\label{eq:U-dist}
 \gamma_k=\left\{ \begin{array}{ll}
   0.0249985, &       k\in \{1, \ldots, 20\}\cup \{81, 100\},\\
   0.000001, & k\in \{21, 22, \ldots, 80\},
                  \end{array}
\right.
\end{equation}
and consider the problem of selecting one of  three best alternatives, i.e., $k=3$.
The optimal value in this problem is
$P(100,3)=0.39711$.
Figure~\ref{fig:Ushape} displays the graphs of sequences
$\{b_{N_{\max}-t-1}\}$ and  $\{U_t(r)\}$, $r=1,2,3$ from which the
form of the stopping region is easily inferred.
\par
Recall that the optimal policy stops when $Y_t=U_t(R_t)> b_{N_{\max}-t-1}$ provided that the decision process
arrives at time $t$. Therefore
the stopping region corresponds to the set of time instances for which
the graphs of $\{U_t(r)\}$, $r=1,2, 3$
are above the graph of $\{b_{N_{\max}-t+1}\}$.
In particular, Figure~\ref{fig:Ushape} shows that the optimal stopping policy is the following.
If
the decision process does not terminate due to horizon randomness then:
pass  the first four observations $t=1,\ldots, 4$;
at time instances  $t=5, \ldots, 15$ stop at the observation with the relative rank one,
if it exists; if not, pass observations $t=16, \ldots, 30$; at time instances $t=31, \ldots, 52$
stop at  the observation with the relative rank one, if it exists; if not,
at time instances $t=53, \ldots, 69$ stop at the observation with the relative rank one
or two, if it exists; if not, at time instances $t=70, \ldots, 99$ stop at
the observation with the relative rank one, two, or three, if it exists; if not, stop at the last observation.

\subsubsection{Expected rank minimization over random horizon}
\label{sec:Pettitt}
In this setting [Problem~(P8) of Section~\ref{sec:selection-problems}]
we would like to minimize the expected absolute rank on the event that the stopping occurs
before $N$; otherwise we receive the absolute rank of the last available observation, $A_{N, N}=R_N$.
Formally, the corresponding stopping problem is
\begin{eqnarray*}
 V_*(N_{\max}) &:=&
 \min_{\tau\in \sT(\sR)} \rE \big[ A_{\tau, N} {\bf 1}\{N\geq \tau\} + R_N {\bf 1} \{N<\tau\}\big]
 \\
& =& - \max_{\tau \in \sT(\sR)} \rE\big [ (R_N-A_{\tau, N}) {\bf 1}\{N\geq \tau\}-R_N]
 \\
 &=& - \max_{\tau \in \sT(\sR)} \rE\big [ (R_N-A_{\tau, N}) {\bf 1}\{N\geq \tau\}] + \frac{1}{2}(1+\rE N).
\end{eqnarray*}
Thus, letting $q(A_{t, k})=R_k-A_{t,k}$ for $t\leq k$ we note that
\begin{eqnarray*}
 I_{t,k}(r)= \rE\big [q(A_{t, k})\,|\,R_1=r_1,\ldots, R_{t-1}=r_{t-1}, R_t=r\big]
 =\frac{1}{2}(k+1) - \frac{k+1}{t+1} r
 \end{eqnarray*}
and therefore
\[
U_t(r)= J_t(r) = \sum_{k=t}^{N_{\max}} \gamma_k I_{t,k}(r) = \Big(\frac{1}{2}-\frac{r}{t+1}\Big)
 \sum_{k=t}^{N_{\max}} (k+1)\gamma_k.
\]
If $N_{\max}=\infty$ then we require that
$\rE N <\infty$; this ensures condition (\ref{eq:u-finite}).
\par
In this setting $\nu=N_{\max}$ or $\nu=\tilde{N}_{\max}$ depending on
support of the distribution of $N$, and
\[
 y_t(j)= \Big(\frac{1}{2}- \frac{j}{t+1}\Big)
 \sum_{k=t}^{\nu}  (k+1) \gamma_k,\;\;\;
 f_t(j)=\frac{1}{t},\;\;\;j=1,\ldots, t,\;\;\;t=1,\ldots, \nu.
\]
The recursion for computation of the optimal value is obtained by substitution of these
formulas in (\ref{eq:alg-recursion}): $b_1=-\infty$,
$b_2= 0$,
and for $t=2,\ldots, \nu$
\begin{eqnarray}
 b_{t+1} &=& \frac{1}{\nu-t+1} \sum_{j=1}^{\nu-t+1}\bigg[b_t \,\vee \,
 \Big(\frac{1}{2} -\frac{j}{\nu-t+2}\Big) \sum_{k=\nu-t+1}^{\nu} (k+1)\gamma_k\bigg].
\nonumber
 \\
&=& b_t + \frac{1}{\nu-t+1}
\sum_{j=1}^{\nu-t+1}\bigg[\Big(\frac{1}{2} -\frac{j}{\nu-t+2}\Big) \sum_{k=\nu-t+1}^{\nu} (k+1)\gamma_k
-b_t\bigg]_+~.
\label{eq:opt-val-expr}
\end{eqnarray}
The optimal policy is to stop at time $t$ if $Y_t=U_t(R_t) > b_{\nu-t+1}$, i.e.,
\begin{align*}
\tau_{\ast}= \min\;
\bigg\{t=1,\ldots,\nu:
\Big(\frac{1}{2}-\frac{R_t}{t+1} \Big)\sum_{k=\nu-t+1}^{\nu}(k+1)\gamma_k> b_{\nu-t+1}
\bigg\}.
\end{align*}
Note that $V_*(N_{\max})=b_{N_{\max}+1} + \frac{1}{2}(1+\rE N)$.
\par
\citeasnoun{Gianini-Pettitt} considered distributions of $N$ with finite right endpoint $N_{\max}$ and
studied asymptotic behavior of
the optimal value $V_*(N_{\max})$ as $N_{\max}\to\infty$. In particular, for distributions satisfying
$\rP(N=k|N\geq k)=(N_{\max}-k+1)^{-\alpha}$, $k=1,\ldots, N_{\max}$, $N_{\max}=1,2,\ldots$ with
$\alpha>0$ one has: (a)~if $\alpha<2$ then $V_*(N_{\max})\to \infty$ as $N_{\max}\to\infty$;
(b)~if $\alpha>2$ then $\lim_{N_{\max}\to\infty} V_*(N_{\max})=3.86945\ldots$;
(c)~if $\alpha=2$ then  $\limsup_{N_{\max}\to\infty} V_*(N_{\max})$ is finite and greater
 than $3.86945\ldots$.
Thus, if $\alpha>2$ then the optimal value $V_*(N_{\max})$ coincides asymptotically with the one
in the classical problem of minimizing the expected rank studied in \citeasnoun{chow}; see Problem~(P4)
in Section~\ref{sec:selection-problems}. On the other hand, if $N$ is uniformly distributed on $\{1, \ldots, N_{\max}\}$,
i.e. $\alpha=1$, then $V_*(N_{\max})\to \infty$ as $N_{\max}\to\infty$.
 \par
We illustrate these results  in Table~\ref{tab:ExpRankUniform}.
The first row of the table, $\alpha=1$, corresponds to the uniform distribution where
$\gamma_k=1/N_{\max}$, $k=1, \ldots, N_{\max}$, while  for general $\alpha>0$
\[
 \gamma_k=\frac{1}{(N_{\max}-k+1)^{\alpha}} \prod_{j=1}^{k-1}
 \bigg[1- \frac{1}{(N_{\max}-j+1)^\alpha}\bigg],\;\;\;k=1,\ldots, N_{\max};
\]
see \citeasnoun{Gianini-Pettitt}.
\begin{table}
\begin{center}
{\small
 \begin{tabular}{||c||c|c|c|c|c|c||}\hline\hline
  $N_{\max}$ & 100 & 500 & $10^3$ & $10^4$ & $10^5$ & $10^6$ \\
\hline
$\alpha=1$ & 4.74437&8.42697 &10.70615 & 23.34298 &50.43062 & 108.71663
\\
$\alpha=2$ & 3.83593 &  4.14133      &  4.18918      & 4.23792 &4.24381 & 4.24444\\
$\alpha=3$&3.61069&3.80588&3.83549& 3.86542 & 3.86909& 3.86947
\\
\hline\hline
 \end{tabular}
}
\caption{ Optimal values $V_*(N_{\max})$ computed using (\ref{eq:opt-val-expr}).
\label{tab:ExpRankUniform}}
\end{center}
\end{table}
It is seen from the table that in the case $\alpha=3$ the optimal value
approaches the universal limit of \citeasnoun{chow} as $N_{\max}$ goes to infinity.
For $\alpha=2$ the formula (\ref{eq:opt-val-expr}) yields the optimal value $4.2444\ldots$; this complements the result of \citeasnoun{Gianini-Pettitt} on boundedness of the optimal value.

\subsection{Multiple choice problems}

The existing literature  treats sequential multiple choice problems
as  problems of multiple stopping.
However, if the reward function has an additive structure, and the involved random variables are independent
then these problems can be reformulated in terms of the sequential assignment problem of
Section~\ref{sec:SAP}.
Under these circumstances
the results of \citeasnoun{DLR} 
are directly applicable
and can be used in order
to construct optimal selection rules.
We illustrate this approach in the next two examples.

\subsubsection{Maximizing the probability of selecting the best observation with $k$ choices}
This setting was first considered by \citeasnoun{GiMo}, and it is discussed in Section~\ref{sec:selection-problems}
as Problem~(P9). The goal is to maximize the probability for selecting the best observation  with $k$ choices, i.e., to maximize
\[
 \rP\big\{\cup_{j=1}^k (A_{\tau_j,n}=1)\big\} = \sum_{j=1}^k \rP(A_{\tau_j,n}=1)
\]
with respect to the stopping times $\tau^{(k)}=(\tau_1,\ldots, \tau_k)$,
$\tau_1<\cdots<\tau_k$ of the filtration $\sR$.
This problem is equivalent to the following version of the sequential assignment problem~(AP1)
[see Section~\ref{sec:SAP}].
\begin{quote}
Let $0=p_1=\cdots=p_{n-k}<p_{n-k+1}=\cdots=p_{n}=1$, and let
\[
 Y_t= \frac{t}{n} {\bf 1}\{R_t=1\}, \;\;t=1, \ldots, n.
\]
The goal is to
maximize  $S(\pi)=\rE \sum_{t=1}^n p_{\pi_t} Y_t$ with respect to $\pi \in \Pi(\sY)$, where
$\Pi(\sY)$ is the set of all non--anticipating policies of filtration $\sY$, i.e.,
\mbox{$\{\pi_t=j\}\in \sY_t$} for all $j=1,\ldots, n$ and $t=1,\ldots, n$.
\end{quote}
The relationship between sequential assignment and multiple choice problems
is evident:
if a policy $\pi$ assigns
$p_{\pi_t}=1$
to the observation $Y_t$ then the corresponding $t$th observation is selected, i.e.,
events $\{p_{\pi_t}=1\}$ and $\cup_{j=1}^k \{\tau_j=t\}$
are equivalent.
\par
The optimal policy for the above assignment  problem is characterized by   Theorem~\ref{th:derman}.
Specifically, for $t=1, \ldots, n$
let $p_{t_1}\leq p_{t_2}\leq \cdots\leq p_{t_{n-t+1}}$
be the subset of the coefficients $\{p_1,\ldots, p_n\}$  that are left unassigned at time $t$.
Let $s_t=\sum_{i=1}^{n-t+1} p_{t_i}$ denote  the
number of observations to be selected
(unassigned coefficients $p$'s equal to $1$).
The optimal policy $\pi_*$ at time $t$ partitions
the real line
by numbers
\[
-\infty=a_{0,n-t+1}\leq a_{1, n-t+1}\leq \cdots\leq a_{n-t, n-t+1}\leq a_{n-t+1, n-t+1}=\infty,
\]
and prescribes to select the $t$th observation   if  $Y_{t}> a_{n-t+1-s_t, n-t+1}$. In words,
the last inequality means that the observation is selected if $Y_t$ is greater than
the $s_t$-th largest number among the numbers
$a_{1, n-t+1}, a_{2, n-t+1}, \ldots, a_{n-t, n-t+1}$.  These numbers are given by the following
formulas: $a_{0, n-t+1}=-\infty$, $a_{n-t+1, n-t+1}=\infty$, and for $j=1, \ldots, n-t$
\[
 a_{j, n-t+1}= \int_{a_{j-1, n-t}}^{a_{j, n-t}} z \rd F_{t+1}(z) + a_{j-1, n-t} F_{t+1}(a_{j-1, n-t})
 + a_{j, n-t} (1-F_{t+1}(a_{j, n-t})),
\]
where $F_t$ is the distribution function of $Y_t$. The optimal value of the problem is
\begin{equation}\label{eq:S-S}
 S_*(k)=S(\pi_*; k)= \sum_{j=1}^k a_{n-j+1, n+1}~.
\end{equation}
\par
In our case $F_t(z)=(1-\frac{1}{t}) {\bf 1}(z\geq 0)+\frac{1}{t} {\bf 1}(z\geq \frac{t}{n})$,
$t=1,\ldots, n$
which yields
\begin{eqnarray}
 a_{j, n-t+1} &=& \tfrac{1}{n} {\bf 1}\big(a_{j-1,n-t}<\tfrac{t+1}{n}\leq a_{j, n-t}\big)
 \nonumber
 \\
&&\;+\; a_{j-1, n-t} \big[\big(1-\tfrac{1}{t+1}\big) {\bf 1}(a_{j-1, n-t}\geq 0)+\tfrac{1}{t+1} {\bf 1}\big(a_{j-1,n-t}\geq \tfrac{t+1}{n}\big)\big]
\label{eq:a-a}
\\
&& \;+\; a_{j, n-t} \big[\big(1-\tfrac{1}{t+1}\big) {\bf 1}(a_{j, n-t} < 0)+\tfrac{1}{t+1}
{\bf 1}\big(a_{j,n-t}< \tfrac{t+1}{n}\big)\big]
\nonumber
\end{eqnarray}
for $j=1, \ldots, n-t$,
$a_{0, n-t+1}=-\infty$, $a_{n-t+1, n-t+1}=\infty$, and by convention we set
\mbox{$-\infty\cdot 0=\infty\cdot 0=0$}.
\par
Table~\ref{tab:a-a} gives optimal values $S_*(k)$ for $n=10^4$ and different $k$. Note that the case $k=1$
corresponds to the classical secretary problem. It is clearly seen that the optimal probability
of selecting the best observtation grows fast  with the number of possible choices~$k$.
The numbers presented in the table agree with those given in Table~4 of
\citeasnoun{GiMo}.
\begin{table}
\begin{center}
{\footnotesize
\begin{tabular}{||c|c|c|c|c|c|c|c|c|c||} \hline\hline
$k$ & 1 &2&3&4&5&6&7&8&25\\ \hline
$S_*(k)$ &0.36791& 0.59106&0.73217 & 0.82319 &0.88263  & 0.92175 &0.94767& 0.96491 & 0.999997\\
\hline\hline
\end{tabular}
}
\caption{Optimal values $S_*(k)$
in the problem of maximizing the probability of selecting the best option with $k$ choices.
The table is  computed using (\ref{eq:a-a}) and (\ref{eq:S-S}) for $n=10^4$.
\label{tab:a-a}}
\end{center}
\end{table}
\par
The structure of the optimal policy allows to compute  distribution of the time required for
the subset selection.
As an  illustration, we consider computation of the expected time required for selecting two
options ($k=2$).
According to the optimal policy the first choice is made
at time
$\tau_1:= \min\{t=1,\ldots, n: Y_t> a_{n-t-1, n-t+1}\}$,
while the second choice  occurs at time
$\tau_2:=\min\{t>\tau_{1}: Y_t> a_{n-t, n-t+1}\}$.
Then the expected time to the subset selection~is
\begin{eqnarray}\label{eq:tau1+tau2}
 \rE \tau_2 &=& \rE \tau_1 + \rE (\tau_2-\tau_1),
\end{eqnarray}
where
\begin{align}
&\rE \tau_1 = 1+ \sum_{j=1}^{n-1} \prod_{t=1}^{j} F_t (a_{n-t-1, n-1+1})
\label{eq:tau-1}
\\
&\rE (\tau_2-\tau_1) =  1+ \sum_{i=1}^{n-2} \rP(\tau_2-\tau_1>i) = 1+
 \sum_{j=1}^{n-1}\sum_{i=1}^{n-j-1} \rP(\tau_2-\tau_1>i\,|\,\tau_1=j) \rP(\tau_1=j)
\nonumber
 \\
&= 1+ \sum_{j=1}^{n-1}\sum_{i=1}^{n-j-1} \prod_{t=1}^{j+i} F_t(a_{n-t, n-1+1}) \rP(\tau_1=j)
\nonumber
\\
&= 1+ \sum_{j=1}^{n-1}\sum_{i=1}^{n-j-1} \prod_{t=j+1}^{j+i} F_t(a_{n-t, n-1+1})
\big[1-F_j(a_{n-j-1, n-j+1})\big] \prod_{t=1}^{j-1} F_t(a_{n-t-1, n-t+1}).
\label{eq:tau-2}
\end{align}
These formulas are clearly computationally amenable and easy to code on a  computer.

\subsubsection{Minimization of the expected average rank with $k$ choices}
In this problem that it is discussed in Section~\ref{sec:selection-problems} as Problem~(P10)
we want to minimize the expected average rank of the $k$ selected observations:
\[
 \min_{\tau^{(k)}} \rE \bigg(\frac{1}{k} \sum_{j=1}^k A_{\tau_j,n}\bigg),
\]
where $\tau^{(k)}=(\tau_1,\ldots, \tau_k)$, $\tau_1<\cdots<\tau_k$ are stopping times of filtration
$\sR$.
\par
This setting is equivalent to the following sequential assignment problem.
\begin{quote}
Let $0=p_1=\cdots=p_{n-k}<p_{n-k+1}=\cdots=p_{n}=1$, and let
\[
 Y_t= -\frac{n+1}{t+1} R_t, \;\;t=1, \ldots, n.
\]
The goal is to
maximize  $S(\pi)=\rE \sum_{t=1}^n p_{\pi_t} Y_t$ with respect to $\pi \in \Pi(\sY)$.
\end{quote}
Note that here
$F_t$ is a discrete distribution with atoms at $y_t(\ell)=-\frac{n+1}{t+1} \ell$, $\ell=1,\ldots, t$
and corresponding probabilities $f_t(\ell):=\rP\{Y_t=y_t(\ell)\}=\frac{1}{t}$.
The structure of the optimal policy is exactly as in the previous section: at time $t$ the real line
is partitioned by
 real numbers $a_{j,n-t+1}$, $j=0,\ldots, n-t+1$ and $t$th option
 if $Y_{t}> a_{n-t+1-s_t, n-t+1}$,
where $s_t$ stands for the number of coefficients $p_i$ equal to $1$ at time $t$.
The
constants
$\{a_{j, n-t+1}\}$ are
determined  by the following formulas:
$a_{0, n-t+1}=-\infty$, $a_{n-t+1, n-t+1}=\infty$, and for $j=2,\ldots, n-t$
\begin{eqnarray*}
  a_{j, n-t+1}&=& \frac{1}{t+1} \sum_{\ell =1}^{t+1} y_{t+1}(\ell)
 {\bf 1}\big\{ y_{t+1}(\ell) \in (a_{j-1, n-t}, a_{j, n-t}]\big\}
 \\
 &&\;+\;
  \frac{a_{j-1, n-t}}{t+1} \sum_{\ell=1}^{t+1} {\bf 1}\big\{y_{t+1}(\ell) \leq a_{j-1, n-t}\big\}
  \;+\; \frac{a_{j, n-t}}{t+1} \sum_{\ell=1}^{t+1} {\bf 1}\big\{y_{t+1}(\ell)  > a_{j, n-t}\big\}~.
\end{eqnarray*}
\par
The optimal value $S_*(k)$ of the problem is again given by (\ref{eq:S-S}).
Table~\ref{tab:b-b} presents $S_*(k)$
for $n=10^5$ and different values of $k$. It worth noting that $k=1$
corresponds to the standard problem of expected rank  minimization [Problem~(P4)] with well known
asymptotics
$S_*(k)\approx 3.8695\ldots$  as $n$ goes to infinity.
Using formulas (\ref{eq:tau1+tau2}), (\ref{eq:tau-1}) and (\ref{eq:tau-2})
we also computed expected time required for $k=2$ selections
when $n=10^3$: $\rE\tau_1\approx 396.25983$ and $\rE \tau_2\approx 610.54822$.  Such performance metrics were not established so far and our approach
illustrates the simplicity with which this can be done.

\begin{table}
\begin{center}
{\footnotesize
\begin{tabular}{||c|c|c|c|c|c|c|c|c|c||} \hline\hline
$k$ & 1 &2&3&4&5&6&7&8&25\\ \hline
$S_*(k)$ &3.86488 &4.50590  &5.12243    &5.72330  &6.31262 & 6.89285 &
7.46574 &8.03255  &17.22753\\
\hline\hline
\end{tabular}
}
\caption{The optimal value $S_*(k)$ in the problem
of minimization of the expected average rank with $k$ choices for $n=10^5$.
\label{tab:b-b}}
\end{center}
\end{table}

\subsection{Miscellaneous problems}
The next two examples illustrate applicability of the proposed framework to some
other problems of optimal stopping.

\subsubsection{Moser's problem with random horizon}
\label{sec:moser}
This is  Problem~(P11) of Section~\ref{sec:selection-problems}.
The stopping problem is
\begin{eqnarray*}
 V_*(N_{\max}) :=
 \max_{\tau\in \sT(\sX)} \rE [ (X_{\tau}- X_N) {\bf 1}\{\tau\leq N\}] + \mu.
\end{eqnarray*}
Define $Y_t= \rE \big[(X_{t}-X_N) {\bf 1}\{t \leq N\} \,|\,\sX_t]$; then
\begin{eqnarray*}
 Y_t= \sum_{k=t}^{N_{\max}} \rE \big[(X_t-X_N) {\bf 1}\{N=k\} \,|\, \sX_t\big]= (X_t-\mu)
 \sum_{k=t}^{N_{\max}} \gamma_k,
\end{eqnarray*}
and for any stopping time $\tau\in \sT(\sX)$
\[
 \rE \big[ (X_\tau - X_N) {\bf 1}\{\tau\leq N\}\big] = \sum_{t=1}^{\infty}
 \rE \Big[{\bf 1}\{\tau=t\} \rE \big\{(X_t-X_N) {\bf 1}\{t\leq N\}\,|\, \sX_t\big\}\Big] = \rE Y_\tau.
\]
Thus, the original stopping problem is equivalent to the problem of stopping the sequence of independent
random variables $Y_t=(X_t-\mu)\sum_{k=t}^{N_{\max}} \gamma_k$, $t=1,\ldots, N_{\max}$, and the optimal
value is
\[
 V_*(N_{\max})= \mu+\max_{\tau\in \sT(\sY)} \rE Y_\tau.
\]
The distribution of $Y_t$ is $F_t(z)=G(\mu+ \frac{z}{\sigma_t})$, $t=1, \ldots, N_{\max}$,
where $\sigma_t:= \sum_{k=t}^{N_{\max}}\gamma_k$.
Then applying Corollary~\ref{cor:opt-stop} we obtain that the optimal stopping rule
is given by
\begin{align*}
&b_1=-\infty,\;\;b_2=\rE Y_{N_{\max}},
\\
& b_{t+1} =\int_{b_t}^\infty z \rd F_{N_{\max}-t+1}(z) + b_t F_{N_{\max}-t+1}(b_t), \;\;\;t=2,\ldots, N_{\max},
 \\
& \tau_*=\min\{1\leq t\leq N_{\max}: Y_t> b_{N_{\max}-t+1}\}.
\end{align*}
\par
In particular, if $G$ is the uniform $[0,1]$ distribution
then
straightforward calculation yields: $b_2=0$ and
\[
 b_{t+1}= \frac{1}{2\sigma_{N_{\max}-t+1}} \Big(b_t +\tfrac{1}{2}\sigma_{N_{\max}-t+1}\Big)^2,\;\;t=2,\ldots,
 N_{\max}.
\]
The optimal value of the problem is $V_*(N_{\max})= b_{N_{\max}+1}+\frac{1}{2}$.
\par
It is worth noting that the case of
$\gamma_k=0$ for all $k=1,\ldots, N_{\max}-1$ and $\gamma_{N_{\max}}=1$
corresponds to the original Moser's problem
with fixed horizon $N_{\max}$. In this case $\sigma_t=1$ for all $t$,
and the above recursive relationship coincides with the one in \citeasnoun{Moser} which is
$E_{t+1}=\frac{1}{2}(1+E_t^2)$ where $E_t=b_t+\frac{1}{2}$.

\subsubsection{Bruss' Odds--Theorem}
This is the  stopping problem (P12) of Section~\ref{sec:selection-problems}. In this setting we have
\begin{equation}
\label{eq:Y-tt}
Y_t:= \rP\{Z_t=1, Z_{t+1}=\cdots=Z_n=0\,|\,\sZ_t\}=
\left\{\begin{array}{ll}
Z_t \prod_{k=t+1}^n q_k, & t=1,\ldots, n-1,\\
Z_t, & t=n,
\end{array}
\right.
\end{equation}
and then
\[
V_*:= \max_{\tau\in \sT(\sZ)} \rP(Z_\tau=1, Z_{\tau+1}=\cdots=Z_n=0\} = \max_{\tau\in \sT(\sY)} \rE Y_\tau.
\]
\par
Thus, the original stopping problem is equivalent to
stopping the sequence $\{Y_t\}$ which is given in (\ref{eq:Y-tt}). Note that
$Y_t$'s are independent, and  $Y_t$ takes two values $\prod_{k=t+1}^n q_k$  and $0$ for $t=1, \ldots, n-1$,
and $1$ and $0$ for $t=n$ with respective
probabilities $p_t$ and $q_t=1-p_t$.
Therefore applying Corollary~\ref{cor:opt-stop}
we obtain that the optimal stopping rule is given by
\begin{equation}\label{eq:Bruss-1}
 \tau_*= \min\bigg\{t=1,\ldots, n: Y_t > b_{n-t+1}\bigg\},
\end{equation}
where
$b_1=-\infty$, $b_2=\rE Y_n = p_n$,
and for $t=2, 3, \ldots,n$
\begin{align}
 \label{eq:Bruss-2}
 b_{t+1}  &= \int_{b_t}^\infty z \rd F_{n-t+1}(z) + b_t F_{n-t+1}(b_t)
= b_t + p_{n-t+1} \bigg[ \prod_{k=n-t+2}^n q_k - b_t\bigg]_+,
\end{align}
where
$[\cdot]_+=\max\{0,\cdot\}$.
The problem optimal value is $V_*=b_{n+1}$.
\par
Now we demonstrate the stopping rule (\ref{eq:Bruss-1})--(\ref{eq:Bruss-2})
is equivalent to the sum--odds--and--stop
algorithm of \citeasnoun{Bruss}.
According to (\ref{eq:Bruss-1}), it is optimal to stop at the first time instance $t\in \{1, \ldots, n-1\}$
such that
$Z_t=1$ and $b_{n-t+1}(\prod_{k=t+1}^n q_k )^{-1} <1$; if such $t$ does not exist then the stopping time is $n$.
Note that
\begin{equation}\label{eq:Bruss-3}
 \frac{b_{n-t+1}}{\prod_{k=t+1}^n q_k}=  \frac{b_{n-t}}{\prod_{k=t+1}^n q_k} +
 \frac{p_{t+1}}{q_{t+1}}\bigg[1 - \frac{b_{n-t}}{\prod_{k=t+2}^n q_k}\bigg]_+,\;\;t=0,1,\ldots, n-2.
\end{equation}
Define $u_s:= b_s (\prod_{k=n-s+2}^n q_k)^{-1}$, $s=2, \ldots, n+1$. It is evident that $\{u_s\}$ is a
monotone increasing sequence, and with this notation
(\ref{eq:Bruss-3}) takes the form
\begin{eqnarray}
 u_{n-t+1} &=& \frac{1}{q_{t+1}} u_{n-t} + \frac{p_{t+1}}{q_{t+1}} (1-u_{n-t})_+,\;\;\;t=0, 1,\ldots, n-2,
 \label{eq:iter}
 \\*[2mm]
 u_2 &=& \frac{p_n}{q_n}.
 \label{eq:iter-1}
\end{eqnarray}
\par
In terms of the sequence $\{u_s\}$ the optimal stopping rule (\ref{eq:Bruss-1}) is the following:
it is optimal to stop at first time $t\in \{1, \ldots, n-1\}$ such that $Z_t=1$ and $u_{n-t+1}<1$; if
such $t$ does not exist then stop at time $n$.
Formally, define
$t_*:=\min\{t=1,\ldots, n-1: u_{n-t+1}<1\}$ if it exists. Then for any $t\in \{t_*, t_*+1. \ldots, n-1\}$
we have $u_{n-t+1}<1$ and
iterating
(\ref{eq:iter})-(\ref{eq:iter-1})
we obtain
\begin{equation}\label{eq:iter-2}
 u_{n-t+1}= u_{n-t}+ \frac{p_{t+1}}{q_{t+1}} = \sum_{k=t+1}^n \frac{p_k}{q_k},\;\;\;t=t_*, t_*+1,\ldots, n-1.
\end{equation}
Therefore (\ref{eq:Bruss-1}) can be rewritten as
\[
 \tau_*= \inf \bigg\{t=1, \ldots, n-1: Z_t=1 \;\;\hbox{and}\;\; \sum_{k=t+1}^n \frac{p_k}{q_k} <1 \bigg\} \wedge n,
\]
where by convention $\inf\{\emptyset\}=\infty$. In order to
compute the optimal value $V_*=b_{n+1}$
we need to determine $u_{n+1}$. For this purpose we note
that the definition of $t_*$ and  (\ref{eq:iter})
imply
\begin{equation}\label{eq:iter-3}
u_{n-t+1} = \frac{u_{n-t}}{q_{t+1}},\;\;\;t=t_*-1, t_*-2, \ldots, 1, 0,
\end{equation}
and,
in view of (\ref{eq:iter-2}),
$u_{n-t_*+1}=\sum_{k=t_*+1}^n (p_k/q_k)$.
Therefore
iterating  (\ref{eq:iter-3}) we have
\[
 u_{n+1}= \bigg(\prod_{j=1}^{t_*} \frac{1}{q_j}\bigg) u_{n-t_*+1}= \bigg(\prod_{j=1}^{t_*} \frac{1}{q_j}\bigg)
 \sum_{k=t_*+1}^n \frac{p_k}{q_k}.
\]
Taking into account that $u_{n+1}= b_{n+1} (\prod_{j=1}^n q_j)^{-1}$ we finally obtain the optimal value of the problem:
\[
 V_*=b_{n+1}= \prod_{j=t_*+1}^n q_j \sum_{k=t_*+1}^n \frac{p_k}{q_k}.
\]
These results coincide with the statement of  Theorem~1 in \citeasnoun{Bruss}.

\section{Concluding remarks}
\label{sec:conc}

We close this paper with several remarks.
\par\medskip
1.~In this paper we show that  numerous problems of sequential selection can be reduced to the
problem of stopping a sequence of independent random variables with carefully specified
distribution functions.
In terms of computational complexity, we cannot assert  that in all cases our approach leads
to a more efficient algorithm than a dynamic programming recursion tailored for a specific problem instance.
However,   in contrast to the latter, in many cases of interest we are able to derive  explicit recursive relationships
that can be  easily implemented; see, e.g., Problem~(P5) that has not been solved to date, or Problems~(P6) and~(P7)
for which our approach provides explicit expressions for computation of optimal policies under  arbitrary distribution of the horizon length.
The conditioning argument
leads to  rules expressed in terms of ``sufficient statistics''; such rules
are very natural, simple, and easy to interpret.
\par\medskip
2.
The proposed framework is applicable to sequential selection problems that can be reduced to settings with
independent observations and additive reward function.
In addition, it is required
that the number of selections to be made is fixed and does not depend on the observations.
As the paper demonstrates, this class is rather broad. In particular, it includes
selection problems with no-information, rank-dependent rewards and fixed or random horizon.
The framework is also applicable to selection  problems with full information when the
random variables $\{X_t\}$ are observable, and the reward for stopping
at time $t$ is a function of the current observation $X_t$ only.
It is worth noting that
in all these problems the optimal policy is of the memoryless threshold
type. In addition, we demonstrate that  multiple choice
problems with fixed and random horizon and additive reward,  as well as
sequential assignment problems with
independent job sizes and random horizon, are also covered by the proposed framework.
In particular, variants of problems (P9), (P10) and (P12) with random horizon
can also be  solved using the proposed approach.
\par\medskip
3.~Although the approach holds for a broad class of sequential selection problems, there
are settings that
do not belong to the indicated class.
For instance, settings
with rank--dependent reward and full information as in \citeasnoun[Section~3]{GiMo} and \citeasnoun{gnedin}
cannot be reduced to optimal stopping of a sequence of independent random variables.
A prominent example of such a setting is
the celebrated Robbins' problem  of
minimizing the expected rank on the basis of full information.
This problem is still open, and  only bounds on the asymptotic optimal value are
available in the literature.
Remarkably,
\citeasnoun{Bruss-Ferguson} show that
no memoryless threshold rule can be optimal in this setting, and
the optimal stopping rule must depend
on the entire  history.
\par\medskip
4. The proposed approach is not applicable to  settings where the number
of selections
is not fixed and depends on the observations.
This class includes problems of maximizing the number of selections subject to
some constraints; for representative publications in
this direction we refer, e.g., to
\citeasnoun{Samuels-Steele}, \citeasnoun{Coffman},
\citeasnoun{Gnedin99}, \citeasnoun{Arlotto15} and references therein.
Another example
is the multiple choice
problem with zero--one reward; see, e.g., \citeasnoun{Rose} and \citeasnoun{Vanderbei} where the problem of
maximizing the probability of selecting  the $k$ best alternatives was considered.
The fact that the results of \citeasnoun{DLR} are not applicable to the latter problem
was already observed by \citeasnoun{Rose} who mentioned this explicitly.

\paragraph{Acknowledgement.}
The authors thank two anonymous reviewers for exceptionally insightful and
helpful reports that led to significant improvements in the paper.
The research was supported by the grants BSF 2010466 and ISF 361/15.

\bibliographystyle{agsm}

\end{document}